\documentclass[10pt]{article}
\usepackage{amsmath,amssymb,amsbsy,amsfonts,amsthm,latexsym,
            amsopn,amstext,amsxtra,euscript,amscd,amsthm}
\usepackage{amssymb}
\usepackage[utf8]{inputenc}
\usepackage[english]{babel}
\usepackage{graphicx}
\usepackage{fancyhdr}
\usepackage{multicol}
\usepackage{multirow}
\usepackage{amsmath, amsthm, amssymb, amsfonts}
\usepackage{float}
\usepackage{enumerate}
\usepackage{cancel}
\usepackage{tikz}
\usepackage{chapterbib}
\usepackage{cite}

\numberwithin{equation}{section}
\usepackage[hidelinks]{hyperref}
\hypersetup{colorlinks=true,linkcolor=blue}

\def\> {\rightarrow}

\def\>{\rightarrow}
\def\S1{\mathbb S^1}

\newtheorem{Definicion}{Definici\'on}[section]

\newtheorem{teo}[Definicion]{Theorem}

\newtheorem{lem}[Definicion]{Lemma}

\newtheorem{obs}[Definicion]{Remark}
\newtheorem{defi}[Definicion]{Definition}

\newcommand{\fin}{\hfill$\blacksquare$}

\usepackage{fancyhdr}
\title{Periodic fractional Ambrosetti-Prodi for one-dimensional problem with drift}
\author{B. Barrios, L. Carrero and A. Quaas}
\date{}
\begin{document}
\maketitle
\abstract{We establish Ambrosetti–Prodi type results for periodic solutions of one-dimensional nonlinear problems with drift term and drift-less whose principal operator is the fractional Laplacian of order $s\in(0,1)$. We establish conditions for the existence and nonexistence of solutions. The proofs of the existence results are based on the sub-supersolution method combined with topological degree type arguments. We also establish a priori bounds in order to get multiplicity results. We also prove that the solutions are $C^{1,\alpha}$ under some regularity assumptions in the nonlinearities, that is, the solutions of equations are classical. We finish the work obtaining existence results for problems with the fractional Laplacian with singular nonlinearity. In particular, we establish an Ambrosetti-Prodi type problem with singular nonlinearities.}

\section{Introduction}
 Motivated by its many applications nonlocal operators have been extensively studied in the last decades.  They appear for instance in continuum mechanics, finance, biology, physics, population dynamics, and game theory among others (see \cite{Laskin, Ralf,fisica2,fisica}).\\
In this work, we deal with the fractional Laplacian operator that can be defined for $s\in (0,1)$ via its multiplier $|\xi|^{2s}$ in Fourier space. For suitable functions like $u\in L^{\infty}(\mathbb{R}^{N})\cap C^{2s+ \beta}_{loc}(\mathbb{R}^{N})$ it can also, be defined by the formula
\begin{equation}\label{fracLapl}
(-\Delta)^s u(x) := {C_{N, s} \mathrm{P.V.} \int_{\mathbb{R}^N} \frac{u(x) - u(z)}{|x - z|^{N + 2s}}dz},
\end{equation}
here $C_{N, s} > 0$ is a well-known normalizing constant and P.V stands for the principal value.\\

Along with all the works, where this operator appears, we found \cite{biswas,ambrosioV,davilaqt} where the fractional Ambrosetti-Prodi type problem is studied, but not in the periodic setting that is the objective of the manuscript at hand. It is fair to mention that not much is known of the existence of periodic solutions in problems that involve the fractional Laplacian (see \cite{quaas,lis,zcg,ambrosioymoli}). In the local case, $s=1$, the seminal work is the one of Ambrosetti and Prodi \cite{ambrosetti}, where the nonlinearity crosses the first Dirichlet eigenvalue. In $1986$, Mawhin in \cite{mawhin} consider  the following equation
\begin{equation*}
  u''(x) + cu'(x) + g(u(x))=t + h(x),
\end{equation*}
where he replaced the crossing condition of the first eigenvalue by
\begin{equation*}
    \lim\limits_{|x|  \to + \infty} g(x)= + \infty,
\end{equation*}
and proved an Ambrosetti-Prodi type result for periodic solutions. Moreover, in $2005$, Mawhin generalized this result for $p$-laplacian in \cite[Theorem 11]{principal}, see also \cite{berger, chia}.\\
In this paper we are going to extend some results proved in \cite{principal} for nonlocal nonlinear equations, like
\begin{equation}\label{1}
   (\triangle) ^{s}u(x) + cu'(x) + g(u(x)) = t + h(x), \qquad\quad x\in \mathbb{R}
\end{equation}
where $(\triangle)^{s}u(x):= - (-\triangle)^{s}u(x)$ and $s\in(0,1)$. Along all the work without loss of generality, we can assume that we are looking for $2\pi$-periodic solutions.\\

As mentioned above, we will focus on the existence of periodic solutions of $(\ref{1})$ for $s\in(0,1)$, which, as far as we know, were not studied until this date. It is important to remark that, one of our most interesting contributions is the case $s\leq 1/2$ since to prove existence we need $C^{1,\alpha}$ regularity where the value of $s$ plays an important role, and for the case, $s=1/2$ is an open problem as far as we know. As it seems natural to prove the existence results we will assume that $h: \mathbb{R} \to \mathbb{R}$  is a $2\pi$-periodic function. Moreover when $c>0$, we will we assume that
\begin{enumerate}[$(H1)$]
  \item $g: \mathbb{R} \to \mathbb{R}$, $h: \mathbb{R} \to \mathbb{R}$ are Lipschitz functions,
\end{enumerate}
 if $s\leq1/2$ and
\begin{enumerate}[$(H2)$]
  \item $g: \mathbb{R} \to \mathbb{R}$, $h: \mathbb{R} \to \mathbb{R}$ are $C^{\alpha}(\mathbb{R})$  functions for some $\alpha\in(0,1)$,
\end{enumerate}
 in the case $s>1/2$.\\
We are also interested in finding periodic solutions to the driftless problem, that is, when supposing that $c=0$. In  this case, when $s> 1/2$ we also assume $(H2)$, and when $s\leq1/2$ we will additionally
\begin{enumerate}[$(H3)$]
  \item $g$ has  a subcritical growth; i.e. also  $$|g(x)| \leq C( 1 + |x|)^{p},$$ where $1< p< 2^{*}_{s} - 1$ with $2^{*}_{s}= \frac{2}{1-2s}$ when $s<1/2$ and $1<p<+\infty$ when $s=1/2$.
\end{enumerate}
We can now state our main result that is given by the following
\begin{teo}\label{princ2}
Assume that
\begin{equation}\label{2}
  \lim\limits_{|x|  \to + \infty} g(x)= + \infty,
\end{equation}
that $h: \mathbb{R} \to \mathbb{R}$ is a $2\pi$-periodic function is satisfied and
\begin{itemize}
  \item $(H1)$ if $s\leq 1/2$ when $c>0$,
  \item $(H2)$ if $s> 1/2$ when $c\geq0$,
  \item $(H2)$-$(H3)$ if $s\leq 1/2$ when $c=0$.
\end{itemize}
Then there exists $t_{1}$ such that the following properties hold
\begin{enumerate}[(i)]
  \item for $t<t_{1}$, $(\ref{1})$ has no $2\pi$-periodic classical solution,
  \item for $t=t_{1}$, $(\ref{1})$ has at least one $2\pi$-periodic classical solution,
  \item for $t>t_{1}$, $(\ref{1})$ has at least two $2\pi$-periodic classical solutions.
\end{enumerate}
\end{teo}

\begin{obs}
  Notice that in the above Theorem $(H2)$ can be replaced by $g$ Lipschitz, but in order not more ask than is necessary we will keep $(H2)$.
\end{obs}

The proof of the Theorem \ref{princ2} is based on the use of the viscosity sub-supersolution method to invert the operator and then the correct use of Leray-Schauder degree type  arguments. We have to bear in mind that we need more regularity of the solutions. We know that determining the regularity of solutions is, in general, a difficult subject to address and depends on the value of $s$. In our case, by using \cite{silvestre} and \cite[Theorem 1.1]{Lsilvestre} we get the $C^{1,\alpha}$ regularity of the solutions in the case $s>1/2$ and $s<1/2$ respectively. The case $s=1/2$ is the most delicate. We did not find a regularity result in the bibliography that would be useful to us, so by using some ideas of \cite{caffarelli} we prove the $C^{1,\alpha}$ regularity directly, see Theorem \ref{regularity} below.\\
We are interested in proving the non-local sub-supersolution method for viscosity solutions, in general, is necessary that $g$ to be monotone, but since we do not have this. We first consider the linear problem and we prove the standard sub-supersolution method, and then for the nonlinear problem, solve the truncated problem using degree arguments, (see Lemma \ref{supsuper}).
In the order to find a second solution, our concern will be obtaining a priori uniform bounds. For the problem with drift term, we find a $L^{2}$ bound for $u'$, this allows us to obtain uniform estimates for $u$. Then, we use Schauder type estimates when $s>1/2$, \cite[Theorem 1.1]{Lsilvestre} when $s<1/2$ and Theorem \ref{regularity} in the case $s=1/2$ and interpolation inequalities in H\"older space for obtaining bounds in $C^{1,\alpha}$ for $s\in(0,1)$. For the driftless problem, we establish $C^{\alpha}$ estimates for the solutions, when $s>1/2$ we consider the problem with solutions of zero mean value, and when $s\leq 1/2$ we use Theorem Liouville after rescaling and compactness, a Gidas Spruck type argument.\\

Our second objective is to combine Ambrosetti- Prodi type problems with nonlinear singularities, that is, we assume that $g$ is unbounded near the origin. As a model, we consider the following equation
\begin{equation}\label{aps1}
  (\triangle)^{s}u(x) + cu'(x) + u(x) + \frac{\beta(x)}{u^{\mu}(x)} = t,
\end{equation}
where $s\in(0,1)$, $c>0$, $\mu \geq 1$, and $\beta$ is Lipschitz continuous and $2\pi$-periodic strictly positive function. In the local case taking $\mu=2$, the equation is an idealized mass-spring model of the electrostatically actuated micro-electro-mechanical system, where  Gutiérrez and Torres in \cite{at} proved the existence and stability of periodic solutions. Related to that our second main theorem is given by the next
\begin{teo}\label{teo2}
Assume $c>0$ and $\beta$ is Lipschitz continuous and $2\pi$-periodic strictly positive function, then there exists $t_{1}>0$ such that $(\ref{aps1})$
so that (i), (ii), and (iii) of the previous Theorem holds.
\end{teo}

It is clear that by the singularity of $g(x, u)= u(x) + \frac{\beta(x)}{u^{\mu}(x)}$ in order to prove the previous result we can not apply Theorem \ref{princ2}, we need to find positive lower bounds in order to obtain a subsolution of the problem. To prove the multiplicity of solutions, we need upper bounds, for this, we find a priori uniform bound a $L^{2}(0,2\pi)$ for $u'$.\\

Organization of the work: in Section $2$, we prove Theorem \ref{princ2}. Section $3$ is devoted to establishing the results of an equation model that combines the Ambrosetti-Prodi problem with singular nonlinearities; that is, we prove Theorem \ref{teo2}. This section also presents the problem with attractive-repulsive singular nonlinearities. Throughout this work, we will  denote by $C>0$ an arbitrary constant that would change from one like to another, and by $$\bar{h}:= 1/2\pi \int\limits_{0}^{2\pi} h(x)dx,$$ the mean value of the function $h$.

\section{The Fractional Ambrosetti-Prodi Problem}
 In this section, we will prove the existence results given in Theorem \ref{princ2}. For this, we prove the existence of viscosity solution, and then, due to that we establish $C^{1,\alpha}$- regularity, we can conclude that the solutions are classical. Therefore, we are interested in the notion of viscosity solutions which requires the punctual evaluation of the equation using appropriate test functions that touch the solution from above or below. From now on we consider $I=(-4\pi,4\pi)$. We have the following
\begin{defi}\label{defivis}
A function $u\in C(\mathbb{R})$ is a viscosity supersolution (resp. subsolution) of $(\ref{1})$ if for every point $x_{0}\in I$, every  $U_{x_0}$ the neighborhood of $x_{0}$ in $I$ and any $\varphi \in C^{2}(\overline{U_{x_{0}}})$ such that
\begin{itemize}
  \item $u(x_{0})=\varphi (x_{0})$,
  \item $\varphi \leq (\text{resp.} \geq) u$ in $U_{x_0}$,
\end{itemize}
then $$(\triangle)^{s}v_{\varphi_{u}}(x_{0}) + cv_{\varphi_{u}}'(x_{0}) + g(v_{\varphi_{u}}(x_{0})) \leq (resp. \geq) t + h(x_{0}),$$ where
\begin{equation*}
  v_{\varphi_{u}}(x):=  \begin{cases}
\varphi(x), & \text{if}\quad x\in U_{x_0},\\
u(x), & \text{if}\quad x\in \mathbb{R}\backslash U_{x_0}.
\end{cases}
\end{equation*}
Finally, $u$ is a viscosity solution of $(\ref{1})$ if it is both a subsolution and a supersolution of $(\ref{1})$ .
\end{defi}

To prove Theorem \ref{princ2} we use the sub-supersolution method for non-proper equations, that is, the following result
\begin{lem}\label{supsuper}
Assume that the assumptions of Theorem \ref{princ2} hold and that there exist $\eta$ and $\beta$, $2\pi$- periodic $C^{1,\alpha}(\mathbb{R})$- functions, for some $\alpha\in(0,1)$, such that $\eta\leq\beta$ and satisfying
\begin{eqnarray*}
(\triangle)^{s}\eta(x) + c\eta'(x) +g(\eta(x)) &\geq& t + h(x), \\
   (\triangle)^{s}\beta(x) + c\beta'(x) + g(\beta(x)) &\leq& t + h(x),
\end{eqnarray*}
in the viscosity sense. Then the equation (\ref{1}) has at least one $2\pi$-periodic viscosity solution $u\in C^{1,\alpha}_{2\pi}(\mathbb{R})$ verifying $\eta(x)\leq u(x) \leq\beta(x)$ for all $x\in \mathbb{R}$.
\end{lem}
 To prove the previous Lemma we need to establish a $C^{1,\alpha}$ regularity results, which will be also useful to find the second solution.\\
Following the ideas of \cite{caffarelli} we prove the $C^{1,\alpha}$ regularity of the solutions in the case $s=1/2$, $c\geq0$ by using the H\"older estimates of the incremental quotients. Indeed we have
 \begin{teo}\label{regularity}
 Let $u$ be a bounded solution of
 \begin{equation*}
    (\triangle)^{1/2}u(x) + cu'(x)= f(x),  \qquad x\in I_{1}:=(-1,1)
 \end{equation*}
 where $f$ is a Lipschitz continuous function. Then there is $0<\alpha<1$ such that $u\in C^{1,\alpha}(-1/2,1/2)$ and
 \begin{equation*}
     \|u\|_{C^{1,\alpha}(-1/2,1/2)} \leq C(\|u\|_{(I_{1})} +  \|f\|_{Lip}(I_{1})),
 \end{equation*}
 for some constant $C>0$, where
 \begin{equation*}
   \|f\|_{Lip}:= \|f\|_{L^{\infty}} + \sup\limits_{0\ne y\in \mathbb{R},x\in \mathbb{R}}\frac{|f(x+y) - f(x)|}{|y|}.
 \end{equation*}
\end{teo}
\begin{dem}
  By \cite[Proposition 2.1]{fernandez} (see also \cite[Theorem 7.2]{Lsilvestreshw}) we have
  \begin{eqnarray}\label{alfa1}
  \|u\|_{C^{\alpha}(I_{1- \delta/2})}   &\leq& C( \|u\|_{C(I_{1})} + \|f\|_{C(I_{1})}),
  \end{eqnarray}
  where $I_{1- \delta/2}= (-(1 - \delta/2), 1 - \delta/2)$ for any $0<\delta<1$. We define now the incremental quotient
  \begin{equation*}
    v_{h}(x):= \frac{u(x + h) - u(x)}{|h|^{1/2}}, \qquad x\in I_{1- \delta/2},
  \end{equation*}
  that is bounded in $I_{1-\delta/2}$ because of $u\in C^{1/2}(I_{1-\delta/2})$. By direct computations, it is clear that $v_{h}$ is a solution of
  \begin{equation*}
    (\triangle)^{\frac{1}{2}}v_{h}(x) + cv_{h}'(x) = \frac{f(x +h) - f(x)}{|h|^{\alpha}}.
  \end{equation*}
  We take now $f_{1}:= \frac{f(x +h) - f(x)}{|h|^{\alpha}}$, which is bounded in $I_{1-\delta/2}$ by using the fact that $f\in C^{\alpha}(I_{1- \delta/2})$. Moreover,
  \begin{equation*}
     \|f_{1}\|_{C(I_{1- \delta/2})} \leq C  \|f\|_{C^{\alpha}(I_{1 - \delta/2})}.
  \end{equation*}
  Applying again \cite[Proposition 2.1]{fernandez} and $(\ref{alfa1})$ we get that
  \begin{eqnarray*}
      \|v_{h}\|_{C^{\alpha}(I_{1 - \delta})} &\leq& C( \|v_{h}\|_{C(I_{1 - \delta/2})} + \|f_{1}\|_{C(I_{1- \delta/2})}), \\
     &\leq& C( \|u\|_{C^{\alpha}(I_{1 - \delta/2})} + \|f\|_{C^{\alpha}(I_{1- \delta/2})}),
  \end{eqnarray*}
   so by \cite[Lemma 5.6]{Luisxavier} we have $u\in C^{\beta}(I_{1 - \delta})$ with $\beta=2\alpha$ and
   \begin{equation*}
       \|u\|_{C^{\beta}(I_{1 - \delta})} \leq  C( \|u\|_{C^{\alpha}(I_{1 - \delta/2})} + \|f\|_{C^{\alpha}(I_{1- \delta/2})}),
   \end{equation*}
   from previous inequality and $(\ref{alfa1})$ we obtain
    \begin{equation*}
       \|u\|_{C^{\beta}(I_{1 - \delta})} \leq  C( \|u\|_{C(I_{1})} +\|f\|_{Lip(I)}).
   \end{equation*}
    Following this iteration in the step $[1/\alpha]$ $u$ is Lipschitz in a small interval, here $[\cdot]$ denotes the integer part function. Thus we can consider the incremental quotient
  \begin{equation*}
    w_{h}(x):= \frac{u(x + h) - u(x)}{|h|},
  \end{equation*}
and by using the same reasoning as before and by choosing the constant $\delta$ appropriately we have $w_{h}\in C^{\alpha}(-1/2,1/2)$ and using the estimates of $u$ and $w_{h}$ obtained by \cite[Proposition 2.1]{fernandez} we get
 \begin{eqnarray*}
      \|w_{h}\|_{C^{\alpha}(-1/2,1/2)} &\leq& C( \|w_{h}\|_{C(I_{1})} + \|f\|_{Lip(I_{1})}),\\
     &\leq& C(\|u\|_{Lip(I_{1})} + \|f\|_{Lip(I_{1})}),\\
     &\leq& C( \|u\|_{C(I_{1})} +\|f\|_{Lip(I_{1})}).
  \end{eqnarray*}
where $C$ is independent of $h$. In particular, this yields
\begin{equation*}
     \|u'\|_{C^{\alpha}(-1/2,1/2)}\leq C( \|u\|_{C(I_{1})} +\|f\|_{Lip(I_{1})}).
\end{equation*}
Therefore, $u\in C^{1,\alpha}(-1/2,1/2)$ and
\begin{eqnarray}\label{12s}
  \|u\|_{C^{1,\alpha}(-1/2,1/2)} &\leq&  \|u\|_{C^{\alpha}(-1/2,1/2)} +  \|u'\|_{C^{\alpha}(-1/2,1/2)} \nonumber \\
  &\leq&C(\|u\|_{C(I_{1})} +\|f\|_{Lip(I_{1})}).
\end{eqnarray}
 as wanted.
\fin
\end{dem}\\

\begin{delem}
Step 1: First of all we consider the following problem
\begin{equation}\label{semilineal}
  (\triangle)^{s}u(x) + cu'(x) - u(x)= f(x),\quad x\in \mathbb{R},
\end{equation}
where $f$ is Lipschitz continuous. This part is standard nowadays but we do it for completeness.

 Define $\eta_{0}=-\|f\|$ and  $\beta_{0}=\|f\|$ that are sub and super solutions of $(\ref{semilineal})$. We define $$A:= \{ z\in C(\mathbb{R}): z \; \text{is a viscosity subsolution of (\ref{semilineal}) with} \; z\leq \beta_{0} \}.$$

Clearly, $A$ is a non-empty set, because of $\eta_{0}\in A$. Let us also consider $$Z(x):=\sup\limits_{z\in A} z(x), \quad x\in I.$$ Our first objective is to prove that $Z$ is a viscosity solution of $\eqref{semilineal}$. We know (see for instance \cite{Nikos}) that $Z$ is a viscosity subsolution of \eqref{semilineal} so, by contradiction, let us suppose that $Z$ is not a supersolution of \eqref{semilineal}. Then there exists $x_{0} \in I$, $U_{x_{0}}$ and $\varphi\in C^{2}(\overline{U_{x_{0}}})$ such that $\varphi\leq Z$ in $U_{x_{0}}$, $\varphi(x_0)= Z(x_{0})$ and
\begin{equation}\label{desigualdad}
   (\triangle)^{s} v_{\varphi_{Z}}(x_{0}) +  cv_{\varphi_{Z}}' (x_{0}) -v_{\varphi_{Z}}(x_{0}) > f(x_{0}).
\end{equation}
Thus there exists $\epsilon_{1}>0$ such that
\begin{equation}\label{desie1}
   (\triangle)^{s} v_{\varphi_{Z}}(x_{0}) +  cv_{\varphi_{Z}}'(x_{0}) -v_{\varphi_{Z}}(x_{0}) - f(x_{0})> \epsilon_{1}.
\end{equation}
Next, we will use test function $\varphi$ to construct a $w\in A$ such that $w$ is above $Z$, if we get our objective, it would  be a contradiction with the definition of $Z$, and this is due to the assumption that $Z$ is not a supersolution of $(\ref{semilineal})$.
We want to prove now that there exists $\delta_{1}>0$  and $\tilde{\epsilon}>0$ such that
\begin{equation}\label{delta1}
\beta_{0}(x) > \varphi(x) +  \tilde{\epsilon},\quad x\in(x_{0} - \delta_{1},x_{0} + \delta_{1})\subseteq U_{x_{0}}.
\end{equation}
In fact, since $\beta_{0} - \varphi$ attains its local minimum at $x_{0}$, there exists $\delta_{1}>0$ such that
\begin{equation}\label{bvs}
  (\beta_{0}- \varphi)(x_{0}) < (\beta_{0}- \varphi)(x), \qquad x\in (x_0 - \delta_{1}, x_0 + \delta_{1}).
\end{equation}
It is easy to prove that $\beta_{0}(x_0)>Z(x_{0})$. Indeed, if $\beta_{0}(x_0)=Z(x_{0})=\varphi(x_{0})$, we have
\begin{equation*}
    (\triangle)^{s} v_{\varphi_{Z}}(x_{0}) +  cv_{\varphi_{Z}}'(x_{0}) - v_{\varphi_{Z}}(x_{0}) \leq f(x_{0}),
\end{equation*}
that contradicts $(\ref{desigualdad})$. Hence by $(\ref{bvs})$,
\begin{equation*}
0 < (\beta_{0}- \varphi)(x), \qquad x\in (x_0 - \delta_{1}, x_0 + \delta_{1}),
\end{equation*}
so that implies that there exists $\tilde{\epsilon}>0$ such that
\begin{equation*}
\tilde{\epsilon} < (\beta_{0}- \varphi)(x), \qquad x\in (x_0 - \delta_{1}, x_0 + \delta_{1}),
\end{equation*}
hence $(\ref{delta1})$ follows.\\
We define now the continuous function
\begin{equation*}
   w(x):= \begin{cases}
Z(x) \quad x\in \mathbb{R}\backslash (x_{0} -\delta_{1}, x_{0} + \delta_{1})\\
\max\{\varphi(x) + \epsilon, Z(x)\},\quad x\in (x_{0} -\delta_{1}, x_{0} + \delta_{1}),
\end{cases}
\end{equation*}
with $\epsilon= \min\{\tilde{\epsilon},\rho,2/3\epsilon_{1}\}$ and $\rho= \min\{(Z- \varphi)(x_{0} - \delta_{1}), (Z- \varphi)(x_{0} + \delta_{1})\}$. It is clear that $w\in A$. In fact, if $w(x)=Z(x)$ then $w\leq \beta_{0}$, and if $w(x)=\varphi(x) + \epsilon$, then $x\in (x_{0} -\delta_{1}, x_{0} + \delta_{1})$ so by $(\ref{delta1})$ we have $w\leq \beta_{0}$. Thus we have to prove that $w$ is a viscosity subsolution. Indeed, let $y_{0}\in I$ and $\phi\in C^{2}(\overline{U_{y_{0}}})$ such that $w\leq \phi$ in $U_{y_{0}}$  and $w(y_{0})=\phi(y_{0})$. We want to prove that $$(\triangle)^{s} v_{\phi_{w}}(y_{0}) +  cv_{\phi_{w}}'(y_{0}) -v_{\phi_{w}}(y_{0}) > f(y_{0}).$$
 If $w(y_{0})=Z(y_{0})=\phi(y_{0})$, then $$0\geq w(x) -\phi(x) \geq Z(x)- \phi(x), \quad x\in U_{y_0},$$ so $$\phi(x) \geq Z(x), \quad x\in U_{y_0},$$ and since $Z$ is a subsolution we have
 \begin{equation*}
  (\triangle)^{s} v_{\phi_{w}}(y_{0}) + cv_{\phi_{w}}'(y_{0}) - v_{\phi_{w}}(y_{0}) \geq (\triangle)^{s} v_{\phi_{Z}}(y_{0}) + cv_{\phi_{Z}}'(y_{0}) - v_{\phi_{Z}}(y_{0}) > f(y_{0}).
 \end{equation*}
Now, if $w(y_{0})= \varphi(y_{0})+ \epsilon=\phi(y_{0})$, then $\phi'(y_{0})= \varphi'(y_{0})$. By definition of $w$, we have $y_{0}\in(x_{0} -\delta_{1}, x_{0} + \delta_{1})$, and
\begin{equation*}
  \varphi(x) + \epsilon \leq w(x) \leq \phi(x), \quad x\in U^{*}:=U_{y_0}\cap(x_{0} - \delta_{1}, x_{0} + \delta_{1}),
\end{equation*}
hence $\varphi + \epsilon - \phi$ attains minimum at $y_{0}$ in $U^{*}$ and since $v_{\varphi_{Z}}\leq v_{\phi_{w}}$ in $\mathbb{R}\backslash U^{*}$ we have
\begin{equation}\label{csub}
   (\triangle)^{s} v_{\phi_{w}}(y_{0}) +  cv_{\phi_{w}}'(y_{0}) - v_{\phi_{w}}(y_{0}) \geq (\triangle)^{s} v_{\varphi_{Z}}(y_{0}) +  cv_{\varphi_{Z}}'(y_{0}) - v_{\varphi_{Z}}(y_{0})- \epsilon,
\end{equation}
since $\varphi\in C^{2}(\overline{U_{x_{0}}})$, then $(\triangle)^{s} v_{\varphi_{Z}}$ is continuous in $U_{x_{0}}$ (see for instance \cite[Lemma $4$.$2$]{caffarelli}), so for $\epsilon/2$ there exists $\delta>0$ such that if $|x_{0} -y_{0}|<\delta$ then
\begin{equation*}
  |(\triangle)^{s} v_{\varphi_{Z}}(y_{0}) + cv_{\varphi_{Z}}'(y_{0}) - v_{\varphi_{Z}}(y_{0})- \epsilon -f(y_{0}) -((\triangle)^{s} v_{\varphi_{Z}}(x_{0})+ cv_{\varphi_{Z}}'(x_{0}) - v_{\varphi_{Z}}(x_{0}) - \epsilon - f(x_{0}))| <\epsilon/2,
\end{equation*}
hence
\begin{equation*}
 (\triangle)^{s} v_{\varphi_{Z}}(y_{0}) +  cv_{\varphi_{Z}}'(y_{0}) - v_{\varphi_{Z}}(y_{0})- \epsilon -f(y_{0}) -((\triangle)^{s} v_{\varphi_{Z}}(x_{0}) + cv_{\varphi_{Z}}'(x_{0}) - v_{\varphi_{Z}}(x_{0}) - \epsilon - f(x_{0})) >-\epsilon/2,
\end{equation*}
and since $\epsilon\leq 2\epsilon_{1}/3$,  together with $(\ref{desie1})$ we have
\begin{eqnarray*}
   (\triangle)^{s} v_{\varphi_{Z}}(y_{0}) + cv_{\varphi_{Z}}'(y_{0}) - v_{\varphi_{Z}}(y_{0}) - \epsilon -f(y_{0}) &>& (\triangle)^{s} v_{\varphi_{Z}}(x_{0}) + cv_{\varphi_{Z}}'(x_{0}) - v_{\varphi_{Z}}(x_{0}) - f(x_{0}) - 3\epsilon/2, \\
   &\geq& \epsilon_{1} - 3\epsilon/2> 0,
\end{eqnarray*}
that is,
\begin{equation*}
  (\triangle)^{s} v_{\varphi_{Z}}(y_{0}) +  cv_{\varphi_{Z}}'(y_{0}) - v_{\varphi_{Z}}(y_{0})- \epsilon - f(y_{0})>0,
\end{equation*}
by $(\ref{csub})$ and the inequality above, we have
\begin{equation*}
    (\triangle)^{s} v_{\phi_{w}}(y_{0}) +  cv_{\phi_{w}}'(y_{0}) - v_{\phi_{w}}(y_{0}) > f(y_{0}),
\end{equation*}
this implies that $w$ is a subsolution of \eqref{semilineal}, so $w\in A$. Thus, as mentioned above, this contradicts the definition of $Z$ and the result follows.\\

Step 2: We prove now the regularity results of \eqref{semilineal}. We consider the cases according to the value of $s$
\begin{itemize}
  \item For $s\in(1/2,1)$, by \cite[Proposition 2.9]{silvestre}, we have that there exists $\alpha\in(0,1)$ that depend of $s$ such that
\begin{eqnarray}\label{rs1}
  \|u\|_{C^{1,\alpha}(\mathbb{R})} &\leq& C (\|u\|_{C(\mathbb{R})} + \|(\triangle)^{s}u\|_{C(\mathbb{R})})\nonumber \\
   &\leq& C(\|u\|_{C(\mathbb{R})}  +  \|cu'\|_{C(\mathbb{R})} +  \|f\|_{C(\mathbb{R})}) \nonumber\\
   &\leq&  C(\|u\|_{C(\mathbb{R})}  + \|u\|_{C^{1}(\mathbb{R})} +\|f\|_{C(\mathbb{R})}),
\end{eqnarray}
where $C=C(|c|, \alpha, s)$. Now, using the Interpolation inequalities \cite[Theorem 3.2.1]{interpol}, we obtain that, for any $\epsilon > 0$, there is a positive constant $C = C(\alpha, |c|,  \epsilon, s)$, such that
\begin{equation}\label{interpol}
  \|u\|_{C^{1}(\mathbb{R})}\leq \epsilon \|u\|_{C^{1, \alpha}(\mathbb{R})} + C\|u\|_{C(\mathbb{R})}.
\end{equation}
Hence, substituting $(\ref{interpol})$ in $(\ref{rs1})$, and choosing for instance $\epsilon=1/2C$ we get that
\begin{equation*}
   \|u\|_{C^{1,\alpha}(\mathbb{R})} \leq  C(\|u\|_{C(\mathbb{R})}+\|f\|_{C(\mathbb{R})}).
\end{equation*}

  \item For $s=1/2$, by Theorem \ref{regularity} we have there exists $\alpha\in(0,1)$ that depend of $s$ such that
  \begin{equation}\label{rs12}
    \|u\|_{C^{1,\alpha}[0,2\pi]} \leq C( \|u\|_{C(\mathbb{R})} + \|f\|_{Lip(I)}).
  \end{equation}

    \item For $s<1/2$, by \cite[Theorem 1.1]{Lsilvestre} we obtain that there exists $\alpha\in(0,1)$ that depend of $s$ such that
\begin{eqnarray*}
   \|u\|_{C^{1,\alpha}(I)}  &\leq&  C(\|u\|_{C(\mathbb{R})} +  \|f\|_{C^{1-2s + \alpha}(I)}).
\end{eqnarray*}
\end{itemize}

Step 3: Solve the truncated equation and conclude. Let $u\in C_{2\pi}(\mathbb{R})$, we define for each $x\in \mathbb{R}$
\begin{equation*}
  F(u(x))= \begin{cases}
-\beta(x) + t + h(x) - g(\beta(x)) , & \text{if $u(x)>\beta(x)$},\\
-u(x)+ t + h(x) - g(u(x)), & \text{if $\eta(x)\leq u(x)\leq \beta(x)$},\\
-\eta(x)+ t + h(x) - g(\eta(x)), & \text{if $u(x)<\eta(x)$.}
\end{cases}
\end{equation*}

First, by Steps $1$ and $2$ we can define a mapping $K:C^{\alpha}_{2\pi}(\mathbb{R}) \to C^{1,\alpha}_{2\pi}(\mathbb{R})$ by $K(z)= u$ where $u$ is a solution of
\begin{equation*}
  (\triangle)^{s}u(x) + cu'(x) - u(x)= z(x).
\end{equation*}
Since for $\beta<\alpha$ the injection of $C^{1,\alpha}_{2\pi}(\mathbb{R})$ into $C^{1,\beta}_{2\pi}(\mathbb{R})$ is compact, consequently $K:C^{\alpha}_{2\pi}(\mathbb{R}) \to C^{1,\beta}_{2\pi}(\mathbb{R})$ is compact. Now, let us define the map
\begin{eqnarray*}
  N: C^{1,\beta}_{2\pi} (\mathbb{R})&\to& C^{\alpha}_{2\pi}(\mathbb{R})  \\
  u &\longmapsto& F(u)
\end{eqnarray*}
Then, we have $K \circ N: C^{1,\beta}_{2\pi} (\mathbb{R})\to C^{1,\beta}_{2\pi}(\mathbb{R})$ is continuous and compact, from now on we denote $K\circ N:= KN$. Notice that $N(C^{1,\beta}_{2\pi} (\mathbb{R}))$ is bounded. Hence, by Schauder's fixed point theorem, $KN$ has a fixed point $u$, i.e., $u$ is a solution of
\begin{equation*}
  (\triangle)^{s}u(x) + cu'(x) - u(x)= F(u(x)) \qquad x\in\mathbb{R}.
\end{equation*}
Now, we prove that $\eta(x)\leq u(x)\leq \beta(x)$ for all $x\in \mathbb{R}$. We only proof that $u(x)\leq \beta(x)$. The other inequality is similar.
We assume, by contradiction, that $\max (u(x) - \beta(x))= u(\bar x) - \beta(\bar x)> 0$, here $\bar x$ is the point where the maximum is attained. So, we have
 $(\triangle)^{s}u(\bar x) - (\triangle)^{s}\beta (\bar x) \leq 0,$ and also
\begin{eqnarray*}
 (\triangle)^{s}u(\bar x) -(\triangle)^{s}\beta(\bar x) &\geq & -cu'(\bar x) + u(\bar x) + F(u(\bar x))+ c\beta'(\bar x) + g(\beta(\bar x)) - t - h(\bar x)\\
 &\geq& -cu'(\bar x) + u(\bar x) -\beta (\bar x) + t + h(\bar x) - g(\beta(\bar x)) +  c\beta'(\bar x) + g(\beta(\bar x)) - t - h(\bar x)\\
 &=& u(\bar x) - \beta(\bar x) > 0.
\end{eqnarray*}
this gives a contradiction. Therefore, $u(x)\leq \beta(x)$.
\fin
\end{delem}\\

The Lemma \ref{supsuper} provides the following result, that is, the existence of the first solution of $(\ref{1})$. By $(\ref{2})$ we can define
\begin{equation*}
  \theta:=\min\limits_{x\in[0,2\pi]z\in\mathbb{R}} (g(z) - h(x))<\infty.
\end{equation*}
\begin{lem}\label{existence1}
Assume that the assumptions of Theorem \ref{princ2} hold, then there exists $t_{1}> \theta$ such that for $t> t_{1}$ the equation $(\ref{1})$ has at least one $2\pi$-periodic solution and for $t< t_{1}$ has no solution.
\end{lem}
\begin{dem}
Let $t^{*}= \max\limits_{x\in[0,2\pi]} (g(0) - h(x))$. Then for $t\geq t^{*}$ we have that $u\equiv0$ is a supersolution for $(\ref{1})$. Now, for $t\geq t^{*}$ from $(\ref{2})$ there exists $R_{t}>0$ such that
\begin{equation*}
  g(u) > t + h(x)\quad \text{whenever}\quad |u|\geq R_{t},\quad x\in\mathbb{R},
\end{equation*}
so that $-R_{t}$ is a subsolution for $(\ref{1})$. Then using Lemma \ref{supsuper} there exists a $2\pi$-periodic solution $u$ of $(\ref{1})$ with $-R_{t}< u(x)<0$. \\
On the other hand, we will use again the method of sup and super solution to proof that if, for $\tilde{t}< t^{*}$ the equation $(\ref{1})$ has a $2\pi$-periodic solution $\tilde{u}$, then it has a $2\pi$-periodic solution for $t\in[\tilde{t}, t^{*}]$. Indeed, we note that $\tilde{u}$ is a super solution for $(\ref{1})$ when $t\in[\tilde{t}, t^{*}]$, since
\begin{equation*}
  (\triangle) ^{s}\tilde{u}(x) + c\tilde{u}'(x) + g(\tilde{u}(x)) =\tilde{ t} + h(x) \leq t + h(x).
\end{equation*}
By $(\ref{2})$ there exists $R_{t}> - \min\limits_{x\in\mathbb{R}} \tilde{u}$ such that
\begin{equation*}
  g(-R_{t}) > t + h(x),
\end{equation*}
so that $-R_{t}< \min\limits_{x\in\mathbb{R}} \tilde{u}$ is a subsolution for $(\ref{1})$. Then, by Lemma \ref{supsuper}, we can conclude that there exists a $2\pi$-periodic solution for $(\ref{1})$. Let us now define
\begin{equation}\label{deft}
  t_{1}=\inf\{ t \in \mathbb{R}: (\ref{1})\; \text{has a} \; 2\pi-\text{periodic solution}\},
\end{equation}
we now want to prove that $t_1\geq\theta$. Suppose by contradiction that if, $u$  is a $2\pi$-periodic solution for $(\ref{1})$ for some $t< \theta$. Then taking $u(x_0) = \min\limits_{x\in\mathbb{R}} u$ and $\varphi\equiv u(x_0)$ since we have that $\varphi'(x_0)=0$ it is follows
\begin{eqnarray*}
   (\triangle) ^{s}v_{\varphi_{u}}(x_0) &=&  t + h(x_0) - g(u(x_0)) \\
   &<& \theta - (g(u(x_0) - h(x_0)) \leq 0,
\end{eqnarray*}
that clearly is a contradiction, because since $\varphi$ is a constant function we have $(\triangle)^{s}v_{\varphi_{u}}(x_0)\geq0$. Therefore $t_{1}\geq \theta$.
\fin
\end{dem}

\subsection{Existence and multiplicity results}

In this section, we prove the existence of a second solution to equation $(\ref{1})$ using the excision of the degree. First of all,
we proof that when $u$ is a possible $2\pi$-periodic solution then a priori uniform bounds hold for $u$.  Reasoning like in step $2$ of proof Lemma \ref{supsuper} we have $u\in C^{1,\alpha}(\mathbb{R})$ for all $s\in (0,1)$. Moreover, for $s>1/2$ by \cite[Proposition 2.8]{silvestre} and gives
\begin{equation*}
  \|u\|_{C^{2s + \alpha}(\mathbb{R})} \leq C( \|u\|_{C(\mathbb{R})} + \|g(u)\|_{C^{\alpha}(\mathbb{R})} + \|h\|_{C^{\alpha}(\mathbb{R})} + t).
\end{equation*}
then $u$ is a classical solution and we have the following bound.

\begin{lem}\label{cotalem}
Assume that the assumptions of Theorem \ref{princ2} hold. Let $t_{2}>t_{1}$, $c>0$ and let $u$ be a $2\pi$-periodic solution of $(\ref{1})$ for $t\in[t_{1},t_{2}]$. Then there exists constant $M=M(t_{2}, h)>0$
such that
\begin{equation}\label{cota}
  \|u\|_{\infty} < M.
\end{equation}
Moreover, there exists $C>0$ such that
\begin{equation*}
  \|u'\|_{L^{2}(0,2\pi)} \leq C \|h\|_{L^{2}(0,2\pi)}.
\end{equation*}
\end{lem}

\begin{dem}
Let $u$ be a $2\pi$-periodic solution for $t\in[t_{1},t_{2}]$. Integrating the equation $(\ref{1})$ over $[0,2\pi]$ we have
\begin{equation*}
  \int\limits_{0}^{2\pi} (\triangle) ^{s}u(x)dx + c\int\limits_{0}^{2\pi} u'(x)dx - \int\limits_{0}^{2\pi} g(u(x))dx = \int\limits_{0}^{2\pi} tdx + \int\limits_{0}^{2\pi} h(x)dx.
\end{equation*}
By \cite[Lemma 2.2]{lis} and using the $2\pi$-periodicity of $u'$ this gives
\begin{equation}\label{igual}
  \frac{1}{2\pi} \int\limits_{0}^{2\pi} g(u(x))dx = t+ \bar{h}.
\end{equation}
Now, multiplying  $(\ref{1})$ by $u'$ and integrating the equation we get
\begin{equation}\label{laecua}
  \int\limits_{0}^{2\pi} (\triangle) ^{s}u(x)u'(x) dx + c \int\limits_{0}^{2\pi} |u'(x)|^{2}dx + \int\limits_{0}^{2\pi} g(u(x))u'(x)dx =  \int\limits_{0}^{2\pi} tu'(x)dx + \int\limits_{0}^{2\pi} h(x)u'(x)dx.
\end{equation}
Since by the periodicity of $u$
\begin{equation*}
 \int\limits_{0}^{2\pi} g(u(x))u'(x)dx = \int\limits_{u(0)}^{u(2\pi)} g(w)dw= 0,
\end{equation*}
and using \cite[Lemma 4.2]{lis},
\begin{equation*}
  \int\limits_{0}^{2\pi} (\triangle) ^{s}u(x)u'(x) dx=0,
\end{equation*}
then
\begin{equation*}
  c \|u'\|_{L^{2}(0,2\pi)}^{2}= \int\limits_{0}^{2\pi} h(x)u'(x)dx.
\end{equation*}
Hence, using a H\"older inequality and the fact that $c>0$ we have
\begin{equation}\label{deri}
   \|u'\|_{L^{2}(0,2\pi)}\leq C  \|h\|_{L^{2}(0,2\pi)}.
\end{equation}
Now, by $(\ref{2})$ there exists $R=R(t_{2}, \bar{h})>0$ such that
\begin{equation*}
  g(u) > t_{2} + \bar{h},
\end{equation*}
whenever $|u|\geq R$. Therefore, if $|u(x)|\geq R$ for all $x\in[0,2\pi]$, we obtain, using $(\ref{igual})$, that
\begin{equation*}
  t + \bar{h}=  \frac{1}{2\pi} \int\limits_{0}^{2\pi} g(u(x))dx > t_{2} + \bar{h}.
\end{equation*}
Thus $t>t_{2}$ which is impossible, so that there exists $x_{0}\in[0, 2\pi]$ such that $|u(x_{0})|< R$. Then by $(\ref{deri})$ we have
\begin{eqnarray*}
  |u(x)| &\leq& |u(x_{0})| + \left|\int\limits_{x_{0}}^{x} u'(z)dz\right|  \\
           &\leq& R+ \sqrt{2\pi} \|u'\|_{L^{2}(0,2\pi)}\\
       &\leq& R + C \|h\|_{L^{2}(0,2\pi)}:= M,
\end{eqnarray*}
as wanted.
\fin
\end{dem}\\

\begin{lem}\label{cotac1al}
Assume that the assumptions of Theorem \ref{princ2} hold. Let $u$ be a solution of $(\ref{1})$ with $t\in [t_{1}, t_{2}]$ and $c>0$, then there exists $N:=N(M,t_{2},h)$ such that
\begin{equation*}
  \|u\|_{C^{1,\alpha}[0,2\pi]}\leq N,
\end{equation*}
where $\alpha$ depends de $s\in(0,1)$.
\end{lem}
\begin{dem}
Let $t\in[t_{1},t_{2}]$ and let $u$ be a $2\pi$-periodic solution of $\ref{1}$. We will divide it by cases:
\begin{itemize}
  \item Let us consider $s\in(1/2,1)$. By \cite[Proposition 2.9]{silvestre} we have that there exists $\alpha\in(0,1)$ that depends the $s$ such that
\begin{eqnarray}\label{cotac1}
  \|u\|_{C^{1,\alpha}(\mathbb{R})} &\leq& C (\|u\|_{C(\mathbb{R})} + \|(\triangle)^{s}u\|_{C(\mathbb{R})})\nonumber \\
   &\leq& C(\|u\|_{C(\mathbb{R})}  +  \|cu'\|_{C(\mathbb{R})} +  \|g(u)\|_{C(\mathbb{R})}+ \|h\|_{C(\mathbb{R})} + t) \nonumber\\
   &\leq&  C(\|u\|_{C(\mathbb{R})}  + \|u\|_{C^{1}(\mathbb{R})} +\|g(u)\|_{C(\mathbb{R})}+ \|h\|_{C(\mathbb{R})} + t),
\end{eqnarray}
where $C=C(|c|, \alpha, s)$. By Lemma \ref{cotalem} there exists $D=D(t_{2},h)>0$ such that
\begin{equation}\label{cotag}
  \|g(u)\|_{C(\mathbb{R})}\leq D,
\end{equation}
and
\begin{equation}\label{cotauc}
  \|u\|_{C(\mathbb{R})}= \|u\|_{C[0,2\pi]}\leq M(t_{2},h).
\end{equation}
Now, using the Interpolation inequalities \cite[Theorem 3.2.1]{interpol}, we obtain that, for any $\epsilon > 0$, there is a positive constant $C = C(\alpha, |c|,\epsilon, s)$, such that
\begin{equation}\label{interc1}
  \|u\|_{C^{1}(\mathbb{R})}\leq \epsilon \|u\|_{C^{1, \alpha}(\mathbb{R})} + C\|u\|_{C(\mathbb{R})}.
\end{equation}
Hence, substituting $(\ref{cotag})$,$(\ref{cotauc})$ and $(\ref{interc1})$ in $(\ref{cotac1})$, using that $t<t_{2}$ and choosing, for instance $\epsilon=1/2C$, we get that
\begin{equation}\label{c1acota}
   \|u\|_{C^{1,\alpha}(\mathbb{R})} \leq C(M+ D + \|h\|_{C(\mathbb{R})} + t_{2}).
\end{equation}
  \item If $s=1/2$, since $g$ is Lipschitz and $u\in C^{1}(I)$ we also have that $g(u)$ is Lipschitz. Then by Theorem \ref{regularity}  we have that there exists $\alpha\in(0,1)$ that depends the $s$ such that
   \begin{equation}\label{1regu}
      \|u\|_{C^{1,\alpha}[0,2\pi]}\leq  C(\|u\|_{C(\mathbb{R})} + \|u'\|_{C(I)}+ \|h\|_{C^{1}(\mathbb{R})}+ t).
   \end{equation}
   By the $2\pi$-periodicity of $u'$ we have $$\|u'\|_{C(I)}=\|u'\|_{C[0,2\pi]},$$ and using again the interpolation inequalities given in \cite[Theorem 3.2.1]{interpol}, we obtain that, for any $\epsilon > 0$, there is a positive constant $C = C(\alpha, |c|,\epsilon, s, \|g\|_{Lip})$, such that
\begin{equation}\label{interc2}
  \|u\|_{C^{1}[0,2\pi]}\leq \epsilon \|u\|_{C^{1, \alpha}[0,2\pi]} + C\|u\|_{C[0,2\pi]},
\end{equation}
   so substituting $(\ref{interc2})$ in $(\ref{1regu})$ together with $(\ref{cotauc})$, $t<t_2$ and choosing for
instance $\epsilon =1/2C$ we obtain
\begin{equation*}
   \|u\|_{C^{1,\alpha}[0,2\pi]} \leq C(M + \|h\|_{C^{1}(\mathbb{R})}+ t_2).
\end{equation*}
\item Let's assume now that $s<1/2$. By using that $u'\in C(I)$ and the fact that $g$ is Lipschitz we clearly get $g(u)\in C^{1-2s + \alpha}(I)$ for some $\alpha\in(0,2s)$.
Hence by $(\ref{cotag})$ if follows
\begin{equation}\label{gc}
  \|g(u)\|_{C^{1-2s + \alpha}(I)}\leq D + C\|u\|_{C^{1}(I)}.
\end{equation}
Since $s<1/2$  by \cite[Teorem 1.1]{Lsilvestre} and $(\ref{gc})$ we obtain that
\begin{eqnarray}\label{ug}
   \|u\|_{C^{1,\alpha([0,2\pi])}}  &\leq&  C(\|u\|_{C(\mathbb{R})} +  \|g(u)\|_{C^{1-2s + \alpha}(I)}+ \|h\|_{C^{1-2s + \alpha}(I)} + t),\nonumber \\
   &\leq&  C(\|u\|_{C(\mathbb{R})} + D + \|u\|_{C^{1}(I)}+ \|h\|_{C^{1-2s + \alpha}(I)} + t_2).
\end{eqnarray}
Using the interpolation inequalities \cite[Theorem 3.2.1]{interpol}, we obtain the inequality $(\ref{interc1})$, and substituting in $(\ref{ug})$ together with $(\ref{cotauc})$ and choosing  $\epsilon=1/2C$ we have
\begin{equation*}
   \|u\|_{C^{1,\alpha}([0,2\pi])} \leq  C(M + D+  \|h\|_{C^{1-2s + \alpha}(I)} + t_2).
\end{equation*}
 \end{itemize}
We can be concluded in all cases that there exists $N=N(M,t_2,h)$ such that
\begin{equation*}
   \|u\|_{C^{1,\alpha}([0,2\pi])} \leq N.
\end{equation*}
\fin
\end{dem}

Let us consider now the equation $(\ref{1})$ when $c=0$, that is,
\begin{equation}\label{csinderi}
   (\triangle) ^{s}u(x) + g(u(x)) = t + h(x), \qquad\quad x\in \mathbb{R}.
\end{equation}

First of all we establish a priori estimates for solutions with zero mean value when $s>1/2$.
\begin{lem}\label{cotazmv}
Assume $(H2)$ and $s>1/2$. For every $C\in\mathbb{R}$, $t_{2}>t_{1}$, with $t_{1}$ was given in Lemma \ref{existence1}, there exists $\tilde{M}:=\tilde{M}(t_{2})$, such that, if $v$ is a bounded  $2\pi$-periodic solution for $t\in[t_{1}, t_{2}]$ of
\begin{equation}\label{esd}
  (\triangle) ^{s}v(x) + g(C+ v(x)) = t + h(x), \; t\in[t_{1},t_{2}]
\end{equation}
 with $\bar{v}=0$, then there exists $\alpha\in(0,1)$ that depends the $s$ such that
\begin{equation*}
  \|v\|_{C^{\alpha}(0,2\pi)}\leq \tilde{M}.
\end{equation*}
\end{lem}

\begin{dem}
Let $t\in[t_{1},t_{2}]$ and $v$ be a bounded solution of $(\ref{esd})$. Integrating the equation over $[0,2\pi]$ and using \cite[Lemma 2.2]{lis} we get
\begin{equation}\label{gh}
  \frac{1}{2\pi}\int\limits_{0}^{2\pi} g(C+ v(x))dx = t + \bar{h}.
\end{equation}
On the other hand, multiplying the equation $(\ref{esd})$ by $v$ and integrating it follows that
\begin{equation}\label{des}
   \int\limits_{0}^{2\pi} v(x)(\triangle) ^{s} v(x) dx + \int\limits_{0}^{2\pi} g(C + v(x))v(x)dx =  \int\limits_{0}^{2\pi} tv(x)dx + \int\limits_{0}^{2\pi} h(x)v(x)dx.
\end{equation}
Since, $(\triangle)^{s}v(x):= - (-\triangle)^{s}v(x)$ and $\bar{v}=0$, we have
\begin{eqnarray*}
 \int\limits_{0}^{2\pi} v(x)(-\triangle) ^{s}v(x) dx  &=&  \int\limits_{0}^{2\pi} g(C + v(x))v(x)dx - \int\limits_{0}^{2\pi} h(x)v(x)dx,\\
                                                   &=& \int\limits_{0}^{2\pi} (g(C + v(x)) - h(x) - \theta) v(x)dx,
\end{eqnarray*}
where $\theta$ was given in Lemma \ref{existence1}. Therefore, by $(\ref{gh})$ and a H\"older inequality, we have
\begin{eqnarray}\label{vs}
  [v]_{H^{s}(0,2\pi)}^{2} &\leq& \|v\|_{\infty}\left(\int\limits_{0}^{2\pi} (g(C + v(x)) - h(x) - \theta)dx\right) \nonumber\\
                      &\leq& \|v\|_{\infty}2\pi(t- \theta)\nonumber\\
                       &\leq& \|v\|_{\infty}2\pi(t_{2}- \theta) \nonumber\\
                      &:=&\tilde{M}(t_{2})\|v\|_{\infty}.
\end{eqnarray}
By Poincare inequality for functions with zero mean value, together with $(\ref{vs})$ it is clear that
\begin{equation}\label{cv1}
  \|v\|_{H^{s}(0,2\pi)}^{2}\leq\tilde{ M}(t_{2})\|v\|_{\infty}.
\end{equation}
Using now that $H^{s}(0,2\pi) \hookrightarrow C^{\alpha}(0,2\pi)$ for $s>1/2$ and some $\alpha\in(0,1)$ we get
\begin{eqnarray*}
   \|v\|_{C^{\alpha}(0,2\pi)}^{2} &\leq&  \tilde{M}(t_{2})\|v\|_{\infty},  \\
                                  &\leq& \tilde{M}(t_{2})\|v\|_{C^{\alpha}(0,2\pi)},
\end{eqnarray*}
obtaining the desirable conclusion.
\fin
\end{dem}

From the previous Lemma we get the last auxiliary result needed to prove Theorem \ref{princ2}; to uniform a priori bounds for the solutions of $(\ref{csinderi})$. Indeed we have the following
\begin{lem}\label{cotasinderi}
Assume that the assumptions of Theorem \ref{princ2} hold. Let $u$ be a solution of $(\ref{csinderi})$ with $t\in [t_{1}, t_{2}]$, then there exists $M:=M(t_{2})$ such that
\begin{equation}\label{cotasind}
  \|u\|_{\infty}\leq M.
\end{equation}
\end{lem}

\begin{dem}
Suppose by contradiction that there exists a sequence $u_{n}$ of solutions of $(\ref{csinderi})$ with $t_{n}\in [t_{1},t_{2}]$ such that $$\|u_{n}\|_{\infty} \to + \infty \; \text{when} \; n \to + \infty.$$

For $s>1/2$, let us write $$u_{n}=\overline{u_{n}} + v_{n},$$ since clearly $\overline{v_{n}}=0$, by Lemma \ref{cotazmv} we have $\|v_{n}\|_{\infty}\leq \tilde{M}$ so that $|\overline{u_{n}}|\to + \infty$ when $n\to +\infty$ and hence $|u_{n}(x)|\to + \infty$ for a.e. $x\in [0,2\pi]$. Thus by $(\ref{2})$ we get
\begin{equation*}
    g(u_{n}(x)) > t_{2} + h(x) \quad  \quad \text{for n large a.e.}\quad x\in [0,2\pi].
\end{equation*}
Therefore, integrating the equation $(\ref{csinderi})$ over $[0,2\pi]$ we obtain
\begin{equation*}
  t_{2} + \bar{h} \geq  t_{n} + \bar{h}= \frac{1}{2\pi} \int\limits_{0}^{2\pi} g(u_{n}(x))dx > t_{2} + \bar{h},
\end{equation*}
giving us a contradiction.\\
For $s\leq1/2$, we want to reduce the problem to positive solutions in order to use some well know Liouville type results. Let $u$ be a periodic solution of $(\ref{csinderi})$ and $x_0\in \mathbb{R}$ a point where the minimum is attained, so we have $(\Delta)^s u(x_0)>0$, consequently
\begin{eqnarray*}
  g(u(x_0)) &=& h(x_0) + t + (-\Delta)^s u(x_0) \\
   &\leq& h(x_0) + t\leq \|h\|_{L^{\infty}} + t:=C.
\end{eqnarray*}
Thus by \eqref{2} there exists $M_0>0$ such that
\begin{equation*}
  -M_0 < u(x_0) <M_0,
\end{equation*}
that implies $u(x)> - M_0$, $x\in\mathbb{R}$. Let define $w:=u + M_0>0 $, that satisfy
\begin{equation}\label{ecuaw}
   (\Delta)^s w(x) + g(w(x) - M_0) = h(x) + t.
\end{equation}
To prove \eqref{cotasind} will be enough to find a priori uniform bounds for $w$. For that suppose, by contradiction, that there exists a sequence $w_{n}$ of solutions of $(\ref{ecuaw})$ with $t_{n}\in [t_{1},t_{2}]$ such that $$\|w_{n}\|_{\infty} \to + \infty \; \text{when} \; n \to + \infty.$$
Let $x_{n}\in\mathbb{R}$ the point where each function $w_{n}$ attains its maximum and we define
$$0<v_{n}(x):= \lambda_{n}w_{n}(\lambda_{n}^{1/\rho} x + x_{n})\quad x\in\mathbb{R},$$ where $\lambda_{n}:= \frac{1}{\|w_{n}\|_{\infty}}$ and $\rho=\frac{2s}{p-1}$. By direct computations, we get that $v_{n}$ is a solution of
\begin{equation*}
   (\triangle)^{s}v_{n}(x) + \lambda_{n}^{1+2s/\rho}g(\lambda_{n}^{-1}(v_{n}(x) - M_{0}\lambda)) = \lambda_{n}^{1 +2s/\rho}(h(\lambda_{n}x + x_{n}) + t_n),
\end{equation*}
where clearly $0 <v_{n}\leq 1$ and $v_{n}(0)=1$. By \cite[Proposition 2.8]{silvestre} there exists a constant $C^{*}>0$ such that
\begin{equation*}
  \|v_{n}\|_{C^{\alpha}(\mathbb{R})}\leq C^{*},
\end{equation*}
for some $\alpha \in (0,1)$ and, therefore, we can extract a convergent subsequence that converges to $v\in C(\mathbb{R})$ such that $0< v\leq 1$ and $v(0)=1$. Moreover by $(H3)$ we also have
\begin{equation*}
   \lim\limits_{n \to + \infty}\frac{g(\lambda_{n}^{-1}(v_{n}(x) - M_0\lambda))}{(\lambda_{n}^{-1}v_{n}(x))^{p}}=C.
\end{equation*}
Hence, since by the election of $\rho$,
\begin{eqnarray*}
  (\triangle)^{s}v_{n}(x) &=& - \lambda_{n}^{1 +2s/\rho}\frac{g(\lambda_{n}^{-1}(v_{n}(x) - M_{0}\lambda))}{(\lambda_{n}^{-1}v_{n}(x))^{p}}(\lambda_{n}^{-1}v_{n}(x))^{p} + \lambda_{n}^{1 + 2s/\rho}(h(\lambda_{n}x + x_{n}) + t_n) \\
                          &=& -\frac{g(\lambda_{n}^{-1}(v_{n}(x) - M_{0}\lambda))}{(\lambda_{n}^{-1}v_{n}(x))^{p}}v_{n}(x)^{p} + \lambda_{n}^{1 + 2s/\rho}(h(\lambda_{n}x + x_{n}) + t_n),
\end{eqnarray*}
taking the limit when $n \to \infty$, since we can assume that $x_{n} \to x_{0}\in \mathbb{R}$,
\begin{equation*}
  (\triangle)^{s}v(x)+ C v(x)^{p}=0.
\end{equation*}
The contradiction follows by \cite[Theorem $4$]{zhuo} when $s<1/2$, and by \cite[Theorem $1$.$2$]{felmer} in the case $s=1/2$.
\fin
\end{dem}\\

\begin{dep}
By Lemma \ref{existence1} we have that $(i)$ is satisfied. We divide the proof of $(ii)$ and $(iii)$ by steps:

Step 1: Let us prove that for $t=t_{1}$, $(\ref{1})$ has at least one $2\pi$-periodic solution when $c>0$.\\
Let $t_{2} > t_{1}$ and  $(t_{n})$ be a sequence in $(t_{1}, t_{2}]$ which converge to $t_{1}$, and let $u_{n}$ be a $2\pi$-periodic solution of $(\ref{1})$ with $t=t_{n}$ given by Lemma \ref{existence1}. In order to pass to the limit, by Lemma \ref{cotac1al}
\begin{equation*}
   \|u_n\|_{C^{1,\alpha}([0,2\pi])} \leq N(M,t_2,h).
\end{equation*}
So by the Arzelá-Ascoli Theorem, we can extract a convergent subsequence in $C^{1}(\mathbb{R})$ such that
\begin{equation*}
  u_n(x) \to u(x), \qquad\qquad u_n'(x) \to u'(x)\quad \text{uniformly},
\end{equation*}
with $u\in C^{1}(\mathbb{R})$. By using \cite[Corollary 4.6]{caffarelli} we conclude that $u$ is a $2\pi$-periodic solution of $(\ref{1})$ for $t=t_{1}$.\\

Step 2: We will obtain the same conclusion as Step $1$ for $c=0$.\\
Let us now consider $t_{2}>t_{1}$ and $\{t_{n}\}$ be a sequence in $(t_{1},t_{2}]$ such that $t_{n} \to t_{1}$ when $n \to +\infty$. We know that for each $t_{n}$ the equation $(\ref{csinderi})$ has a bounded solution $u_{n}$. Since by Lemma \ref{cotasinderi}  there exists $\tilde{D}=\tilde{D}(t_2)$ such that $$ \|g(u_n)\|_{C(\mathbb{R})}\leq \tilde{D},$$ and \cite[Proposition 2.9]{silvestre} we have that there exists $\alpha\in(0,1)$ that depends the $s$ such that
\begin{eqnarray*}
  \|u_n\|_{C^{\alpha}(\mathbb{R})} &\leq& C (\|u_n\|_{C(\mathbb{R})} + \|(\triangle)^{s}u_n\|_{C(\mathbb{R})})\nonumber \\
   &\leq& C(\|u_n\|_{C(\mathbb{R})} +  \|g(u_n)\|_{C(\mathbb{R})}+ \|h\|_{C(\mathbb{R})} + t_n) \nonumber\\
   &\leq&  C(M(t_2)  + \tilde{D} + \|h\|_{C(\mathbb{R})} + t_2),
\end{eqnarray*}
 up to subsequence, $u_{n}$ converges to some $u\in C^{\beta}(0,2\pi)$ with $0< \beta < \alpha$. From  \cite[Corollary 4.6]{caffarelli} we conclude that $u$ is a bounded $2\pi$-periodic solution of $(\ref{csinderi})$ for $t=t_{1}$.\\

Step 3: Let us prove for $c>0$, the existence of at least two $2\pi$-periodic solutions. For that, let $t_{2}>t_{1}$, and let $\tilde{t}\in(t_{1},t_{2})$ fixed. We can define a map $K:C^{\alpha}_{2\pi}(\mathbb{R}) \to C^{1,\alpha}_{2\pi}(\mathbb{R})$ by $K(z):= u$ where $u$ is solution of
\begin{equation*}
   (\triangle) ^{s}u(x) + cu'(x) = z(x).
\end{equation*}
 Since for $\tilde{\alpha}<\alpha$ the injection of $ C^{1,\alpha}_{2\pi} (\mathbb{R})$ into  $C^{1,\tilde{\alpha}}_{2\pi} (\mathbb{R})$ is compact, consequently $K: C^{\alpha}_{2\pi}(\mathbb{R})\to C^{1,\tilde{\alpha}}_{2\pi} (\mathbb{R})$ is compact.\\  Let us define now $N: C^{1}_{2\pi}(\mathbb{R}) \to  C^{\alpha}_{2\pi} (\mathbb{R})$ by $Nu:= t + h(\cdot) - g(u(\cdot))$. It is clear that, $KN : C^{1}_{2\pi} (\mathbb{R})\to C^{1,\tilde{\alpha}}_{2\pi} (\mathbb{R})$ is also compact.
 Then, since by $(i)$ $(\ref{1})$ has no $2\pi$-periodic solution for $t<t_{1}$, we get, for all $t\leq t_2$
\begin{equation}\label{d1}
  deg( I - KN, \Omega, 0)= 0.
\end{equation}
 $$\Omega= \{ u\in C^{1}_{2\pi}(\mathbb{R}): \|u\|_{C(\mathbb{R})}< M(t_2,h),\quad\|u'\|_{C(\mathbb{R})}< N(M,t_2,h)\}.$$ with $M(t_2,h)$ and $N(M,t_2, h)$ is the bound given by Lemma \ref{cotalem} and \ref{cotac1al} respectively.
On the other hand, by part $(i)$ we consider $u_{1}$ a solution of $(\ref{1})$ for $t_{1}$ that, since
\begin{equation*}
     (\triangle) ^{s}u_{1}(x) + cu_{1}'(x) - g(u_{1}(x)) = t_{1} + h(x)< \tilde{t} + h(x),
\end{equation*}
is a strict $2\pi$-periodic supersolution of $(\ref{1})$ for $\tilde{t}$.
Reasoning like in the proof of Lemma \ref{existence1} we get that there exists $R_{t_{2}}$ such that
\begin{equation*}
  g(-R_{t_2})>t_2 + h(x) >\tilde{t} + h(x), \qquad x\in \mathbb{R}.
\end{equation*}
So that, $-R_{t_{2}}$ is a strict $2\pi$-periodic subsolution for $(\ref{1})$ with $\tilde{t}$. Define
$\tilde{N}: C^{1}_{2\pi} (\mathbb{R})\times [0,1] \to  C^{\alpha}_{2\pi} (\mathbb{R})$ by
$$\tilde{N}(u,\lambda):= \lambda t + \lambda h(\cdot) - \lambda g(u(\cdot)).$$
From the invariance of the degree concerning a homotopy for $t=\tilde{t}$ we have
\begin{equation}\label{deghop}
    deg(I- K\tilde{N}(\cdot,1), \Omega_{1}, 0)=   deg(I- K\tilde{N}(\cdot,0), \Omega_{1}, 0)=   deg(I, \Omega_{1}, 0)=1
\end{equation}
where $$\Omega_{1}:= \{ u \in C^{1}_{2\pi}(\mathbb{R}): - R_{t_{2}} < u(x) < u_{1}(x), \;  x\in\mathbb{R}\quad \|u'\|_{C(\mathbb{R})}< N(M,t_2,h)\},$$
then for $t=\tilde{t}$ we get
\begin{equation}\label{deg1h}
   deg(I- K\tilde{N}(\cdot,1), \Omega_{1}, 0)= deg(I- KN, \Omega_{1}, 0)= 1.
\end{equation}
 Since we can take $M(t_{2},h)$ large enough such that $\Omega_{1} \subseteq \Omega$, by using $(\ref{d1})$, $(\ref{deg1h})$ and the property excision of the Leray-Schauder degree we have
\begin{equation*}
  deg( I -KN, \Omega\backslash \overline{\Omega_{1}}, 0) = deg( I-KN, \Omega, 0) - deg( I- KN, \Omega_{1}, 0) =-1,
\end{equation*}
that implies the existence of a $2\pi$-periodic solution of equation $(\ref{1})$ contained in $\Omega\backslash \overline{\Omega_{1}}$ for $\tilde{t}$. Since, $t_{2}> t_{1}$ and $\tilde{t}\in [t_{1}, t_{2}]$ are arbitrary, the proof is complete.\\

Step 4: Finally we obtain the conclusion of $(iii)$ the driftless problem. For that, we take $\tilde{t}\in (t_{1}, t_{2})$ fixed but arbitrary. By \cite[Lemma 3.1]{quaas} we can define $K:C^{\alpha}_{2\pi}(\mathbb{R})\to C^{2s+\alpha}_{2\pi}(\mathbb{R})$ by $K(z):=u$  where $u$ is solution of
\begin{equation*}
  (\triangle)^{s}u(x) + u(x) = z(x).
\end{equation*}
Because the injection of $C^{2s + \alpha}$ into $C^{\alpha}$ is compact, we also get that $K: C_{2\pi}^{\alpha}(\mathbb{R}) \to  C^{\alpha}_{2\pi} (\mathbb{R})$ is compact. By step $2$, we can conclude that there exists $N(M,h,t_2)$ such that
\begin{equation*}
  \|u\|_{C^{\alpha}_{2\pi} (\mathbb{R})} \leq N.
\end{equation*}
 Let us define now $$\Omega:= \{  u\in  C^{\alpha}_{2\pi} :  \|u\|_{C_{2\pi}}< M(t_{2}),\quad [u]_{\alpha}< N\},$$ and $N:C^{\alpha}_{2\pi}  \to  C^{\alpha}_{2\pi} (\mathbb{R})$ by
\begin{equation*}
  Nu= u(\cdot) + t + h(\cdot) - g(u(\cdot)).
\end{equation*}
It is clear that, $KN: C^{\alpha}_{2\pi}  \to  C^{\alpha}_{2\pi} (\mathbb{R})$ is also compact. Let now define Since by $(i)$ the equation $(\ref{csinderi})$ has no solution for $t<t_{1}$ we also have,
\begin{equation}\label{q1}
  deg( I - KN, \Omega, 0)= 0,
\end{equation}
for every $t\leq t_{2}$. On the other hand, by using the method of sup and supersolutions we find a solution in a subset of $\Omega$. In fact, by $(ii)$ is clear that $u_{1}$ is a solution of $(\ref{csinderi})$  for $t_{1}$, so is a strict $2\pi$-periodic supersolution of $(\ref{csinderi})$ for $\tilde{t}$. Doing similarly as in the  proof of Lemma \ref{existence1} we know that, there exists $R_{t_{2}}$ such that
\begin{equation*}
  g(-R_{t_{2}})> t_{2} + h(x)> \tilde{t} + h(x), \qquad x\in \mathbb{R},
\end{equation*}
that also implies that, $-R_{t_{2}}$ is a strict $2\pi$-periodic subsolution for $(\ref{csinderi})$ for $\tilde{t}$. Therefore, reasoning like in the previous step for $t=\tilde{t}$ we get
\begin{equation}\label{deg2}
  deg(I- KN, \Omega_{1}, 0)= 1,
\end{equation}
where $$\Omega_{1}= \{ u \in C^{\alpha}_{2\pi}(\mathbb{R}):  - R_{t_{2}} < u(x) < u_{1}(x),\;  x\in\mathbb{R} \quad [u]_{\alpha}< N \}.$$ As in the previous by step, by taking $M(t_{2})$ large enough that guarantee that $\Omega_1 \subseteq \Omega$, and by using the excision property, we conclude for $t=\tilde{t}$ that
\begin{equation*}
  deg( I -KN, \Omega\backslash \overline{\Omega_{1}}, 0) = deg( I-KN, \Omega, 0) - deg( I- KN, \Omega_{1}, 0) =-1,
\end{equation*}
that implies the existence of a $2\pi$-periodic solution of equation $(\ref{csinderi})$ contained in $\Omega\backslash\overline{\Omega_{1}}$ for $\tilde{t}$. Again, since $\tilde{t}$ and $t_{2}$ are arbitrary, the proof is complete.\\
In the previous steps, we obtained the existence of viscosity solutions, and we know (for instance Lemma \ref{supsuper}) that these solutions with $c\geq0$ are $C^{1,\alpha}(\mathbb{R})$. This implies for the case $s\leq1/2$ that the solutions of the equation are classical. Moreover, for $s>1/2$ by \cite[Proposition 2.8]{silvestre} and $(\ref{c1acota})$ gives
\begin{equation*}
  \|u\|_{C^{2s + \alpha}(\mathbb{R})} \leq C( M + D + \|h\|_{C^{\alpha}(\mathbb{R})} + t_2).
\end{equation*}
Therefore, the solutions of the equation are classical for $s\in(0,1)$.
\fin
\end{dep}

\begin{obs}
\begin{enumerate}[1.]
  \item An alternative proof of the part $(ii)$ and $(iii)$ for the case $s>1/2$ and $c=0$ is to solve the problem with solutions of zero mean value, that is, we consider the following problem
\begin{equation}\label{ecuavm}
  (\triangle) ^{s}v(x) + g(C+ v(x)) = t + h(x),
\end{equation}
where $\bar{v}=0$ and $C$ is fixed but arbitrary. For this, reasoning like in the previous proof and using the Lemma \ref{cotazmv} the result is concluded.
  \item It is fair to mention that in the driftless problem, we can take $\Omega$ the same form as in step $3$, but since it is not necessary that much regularity, we worked in $C^{\alpha}$.
\end{enumerate}
\end{obs}
\section{Ambrosetti-Prodi problem with singular nonlinearities}
  We prove now the existence of positive solutions related to the equation $(\ref{aps1})$; that is Theorem \ref{teo2}. For that we take $$g(x, u):= u(x) + \frac{\beta(x)}{u^{\mu}(x)}.$$ It is good to be in mind that we can not apply Theorem \ref{princ2}, due to the singularity of $g$.\\
To prove the result of existence we know that it is essential to find a lower bound, for this purpose we follow the ideas of \cite{at}.
\begin{lem}\label{infcota}
Let $t>0$ and $u$ be a solution of $(\ref{aps1})$, then there exists $r_t=r_t(t)>0$ such that $u(x)>r_t$ for every $x\in [0,2\pi]$.
\end{lem}
\begin{dem}
Let $x_0$ be the point where the minimum of $u$ in $[0,2\pi]$ is attained. Then
\begin{equation*}
   u(x_0) + \frac{\beta(x_0)}{u^{\mu}(x_0)}\leq t,
\end{equation*}
that in particular implies
\begin{equation*}
 \left(\frac{\beta_{min}}{t}\right)^{1/\mu} \leq u(x_0),
\end{equation*}
where $\beta_{min}:= \min\limits_{x\in[0,2\pi]} \beta(x)$. Taking $r_t:= (\frac{\beta_{min}}{t})^{1/\mu}$, the proof is complete.
\fin
\end{dem}\\
\begin{lem}\label{supaps}
 Assume $t>0$, $c>0$ and let $u$ be a $2\pi$-periodic solution of $(\ref{aps1})$. Then there exists constant $M=M(t)>0$
such that
\begin{equation}\label{cotAs}
  \|u\|_{\infty} < M.
\end{equation}
Moreover, there exists $C>0$ such that
\begin{equation*}
  \|u'\|_{L^{2}(0,2\pi)} \leq C(\|\beta\|_{L^{2}(0,2\pi)}  + t\sqrt{2\pi}).
\end{equation*}
\end{lem}
\begin{dem}
Let $u$ be a $2\pi$-periodic solution for $t>0$. Integrating the equation $(\ref{aps1})$ over $[0,2\pi]$ we have
\begin{equation*}
  \int\limits_{0}^{2\pi} (\triangle) ^{s}u(x)dx + c\int\limits_{0}^{2\pi} u'(x)dx + \int\limits_{0}^{2\pi} u(x) dx + \int\limits_{0}^{2\pi} \frac{\beta(x)}{u^{\mu}(x)}dx= \int\limits_{0}^{2\pi} tdx.
\end{equation*}
By \cite[Lemma 2.2]{lis} and using the $2\pi$-periodicity of $u'$ this gives
\begin{equation}\label{crt}
  \frac{1}{2\pi}\left(\int\limits_{0}^{2\pi} u(x) dx + \int\limits_{0}^{2\pi} \frac{\beta(x)}{u^{\mu}(x)}dx\right)= t.
\end{equation}
Now, there exists $R_t=R(t)>0$ such that
\begin{equation*}
   g(x,u)> t,
\end{equation*}
whenever $u\geq R_t$. Therefore, if $u(x)\geq R_t$ for all $x\in[0,2\pi]$, we obtain, using $(\ref{crt})$, that
\begin{equation*}
  t =  \frac{1}{2\pi} \int\limits_{0}^{2\pi} u(x) dx + \int\limits_{0}^{2\pi} \frac{\beta(x)}{u^{\mu}(x)}dx > t.
\end{equation*}
Thus $t>t$ which is impossible, so that there exists $x_{1}\in[0, 2\pi]$ such that $u(x_{1})< R_t$. Then we have
\begin{eqnarray}\label{cotud}
  u(x) &\leq& u(x_{1}) + \int\limits_{x_{1}}^{x} |u'(z)|dz  \\
           &\leq& R_t+ \sqrt{2\pi} \|u'\|_{L^{2}(0,2\pi)}.
\end{eqnarray}
Let now $x_0$ be the point where the minimum of $u$ in $[0,2\pi]$ is attained. Multiplying $(\ref{aps1})$ for $u'$ and integrating the equation
\begin{equation*}
    \int\limits_{0}^{2\pi}(\Delta)^s u(x)u'(x)dx + c \int\limits_{0}^{2\pi}|u'(x)|^{2}dx + \int\limits_{0}^{2\pi} u(x)u'(x) dx  + \int\limits_{0}^{2\pi}\frac{\beta(x)}{u^{\mu}(x)}u'(x)dx = \int\limits_{0}^{2\pi} tu'(x)dx.
  \end{equation*}
  By \cite[Lemma 4.2]{lis} and using H$\ddot{o}$lder inequality together with the Lemma \ref{infcota}  gives
  \begin{eqnarray*}
     c\|u'\|_{L^{2}(0,2\pi)}^{2} &\leq&\int\limits_{0}^{2\pi}\frac{|\beta|}{u^{\mu}(x_{0})}|u'| + \int\limits_{0}^{2\pi} t |u'| \\
     &\leq& C\|u'\|_{L^{2}(0,2\pi)} (\|\beta\|_{L^{2}(0,2\pi)}  + t\sqrt{2\pi}),
  \end{eqnarray*}
  so that,
\begin{equation*}
  \|u'\|_{L^{2}(0,2\pi)}\leq C(\|\beta\|_{L^{2}(0,2\pi)}  + t\sqrt{2\pi}),
\end{equation*}
  which, together with $(\ref{cotud})$ implies
  \begin{equation*}
    u(x) \leq C(R_t+ \|\beta\|_{L^{2}(0,2\pi)}  + t)=: M
  \end{equation*}
\fin
\end{dem}

In the following result, we prove the existence of the first solution. The proof differs from Lemma \ref{existence1} in the form that the subsolution is obtained and this is due to the singularity of $g$. Define
\begin{equation*}
  \theta:= (\mu \beta_{min})^{1/(\mu +1)}, \qquad t^{*}:= (\mu +1)\left( \dfrac{\beta_{max}}{\mu^{\mu}}\right)^{1/(\mu +1)},
\end{equation*}
where $\beta_{max}:= \max\limits_{x\in[0,2\pi]} \beta(x)$.
\begin{lem}\label{nosol}
There exists $t_{1}> \theta$ such that for $t> t_{1}$ the equation $(\ref{aps1})$ has at least one $2\pi$-periodic solution and for $t< t_{1}$ has no solution.
\end{lem}
\begin{dem}
We first proof that the equation $(\ref{aps1})$ has at least one $2\pi$-periodic solution for $t> t^{*}$. For this, we known by Lemma \ref{supsuper} it suffices find sub and supersolutions $\eta_0$ and $\beta_0$ such that $\eta_0\leq \beta_0$. Therefore, for $t> t^{*}$ we have $\beta_0\equiv (\mu\beta_{max})^{1/(\mu +1)}$ is supersolution. Indeed,
\begin{equation*}
   \beta_0 + \frac{\beta(x)}{\beta_0^{\mu}} \leq \beta_0 + \frac{\beta_{max}}{\beta_0^{\mu}}= t^{*}\leq t.
\end{equation*}
Now, for $t> t^{*}$ we have that $\eta_0\equiv r_t$ is a subsolution. Indeed,
\begin{equation*}
   r_t + \frac{\beta(x)}{r_t^{\mu}} \geq r_t + \frac{\beta_{min}}{r_t^{\mu}}= r_t + t\geq t.
\end{equation*}
It is easily checked that $\eta_0\leq \beta_0$. So, there exists a $2\pi$-periodic solution $u$ of $(\ref{aps1})$ with $\eta_0< u<\beta_0$.
Defining,
\begin{equation*}
   t_{1}:=\inf\{ t \in \mathbb{R}: (\ref{aps1})\; \text{has a} \; 2\pi-\text{periodic solution}\},
\end{equation*}
we will proof that $(\ref{aps1})$ has no $2\pi$-periodic solution for $t<\theta$. Indeed, if $u$ is a $2\pi$-periodic solution for $(\ref{aps1})$ for $t>t_1$, then integrating the equation $(\ref{aps1})$ over $[0,2\pi]$ and using the $2\pi$-periodicity of $u'$ together with  \cite[Lemma 2.2]{lis} we have
\begin{equation}\label{crt}
  \frac{1}{2\pi}\left(\int\limits_{0}^{2\pi} u(x) dx + \int\limits_{0}^{2\pi} \frac{\beta(x)}{u^{\mu}(x)}dx\right)= t.
\end{equation}
Note that the function $$g_{0}(u)=u + \frac{\beta_{min}}{u^{\mu}}$$ reaches its global minimum at $\theta$, so that
\begin{equation*}
  \theta \leq \frac{1}{2\pi}\left(\int\limits_{0}^{2\pi}u(x) dx +\frac{\beta_{min}}{u^{\mu}(x)}dx\right)\leq \frac{1}{2\pi}\left(\int\limits_{0}^{2\pi}u(x) dx +\frac{\beta(x)}{u^{\mu}(x)}dx\right) =t.
\end{equation*}
Therefore, $t_1\geq \theta$.
\fin
\end{dem}\\

\begin{desing}
Part $(i)$ follows from Lemma \ref{nosol}. Now the part $(ii)$, let us consider  $t_{2} > t_{1}$ and  $\{t_{n}\}_{n\in \mathbb{N}}$ to be a sequence in $(t_{1}, t_{2}]$ which converge to $t_{1}$, and let $u_{n}$ be a $2\pi$-periodic solution of $(\ref{aps1})$ with $t=t_{n}$ whose existence is guaranteed by part $(i)$. By Lemma \ref{infcota} and Lemma \ref{supaps} there exists $D=D(t_{2}, \beta)>0$ such that
\begin{equation}
  \|g(u)\|_{C(\mathbb{R})}\leq D,
\end{equation}
  and since $g$ is Lipschitz, we can reason similarly to the proof of Lemma \ref{cotac1al} to obtain there exists $\alpha\in(0,1)$ such that
  \begin{equation*}
   \|u_n\|_{C^{1,\alpha}([0,2\pi])} \leq  N,
\end{equation*}
with $N=N(R_{t_2},D,\beta_{max},t_2)$. By the Arzelá-Ascoli Theorem, we can extract a convergent subsequence in $C^{1}(\mathbb{R})$ such that
\begin{equation*}
  u_n(x) \to u(x), \qquad\qquad u_n'(x) \to u'(x)\quad \text{uniformly},
\end{equation*}
with $u\in C^{1}(\mathbb{R})$. By using \cite[Corollary 4.6]{caffarelli} we conclude that $u$ is a $2\pi$-periodic solution of $(\ref{aps1})$ for $t=t_{1}$.\\
 Finally, the part $(iii)$, let $t_2>t_1$ and let $\tilde{t}\in (t_1,t_2)$ fixed. We can define a map $K: L^{\infty}(\mathbb{R}) \to C^{1,\alpha}_{2\pi} (\mathbb{R})$ by $K(z):=u$ where $u$ is solution of
 \begin{equation*}
      (\triangle) ^{s}u(x) + cu'(x) =z(x).
 \end{equation*}
 Since for $\tilde{\alpha}<\alpha$ the injection of $ C^{1,\alpha}_{2\pi} (\mathbb{R})$ into  $C^{1,\tilde{\alpha}}_{2\pi} (\mathbb{R})$ is compact, consequently $K: L^{\infty}(\mathbb{R}) \to C^{1,\tilde{\alpha}}_{2\pi} (\mathbb{R})$ is compact. Now, let us define the map
 \begin{eqnarray*}
  N: C^{1}_{2\pi} (\mathbb{R}) &\to&  L^{\infty}(\mathbb{R})  \\
  u &\longmapsto& t - u - \frac{\beta}{u^{\mu}}.
\end{eqnarray*}
Then $KN:  C^{1}_{2\pi} (\mathbb{R}) \to C^{1,\tilde{\alpha}}_{2\pi} (\mathbb{R})$ is continuous and compact. Since by $(i)$ $(\ref{aps1})$ has no $2\pi$-periodic solution for $t<t_1$, we also have,
\begin{equation}\label{deg0}
  deg( I - KN, \Omega, 0)= 0,
\end{equation}
for every $t\leq t_{2}$, where $\Omega=\{u\in C^{1}_{2\pi} (\mathbb{R}):  r_{t_2} < u(x) < R_{t_2},\;  x\in \mathbb{R}\quad   \|u'\|_{\infty}<N\}.$  On the other hand, by $(ii)$ we consider $u_{1}$ is a solution of $(\ref{aps1})$ for $t_{1}$ so that
\begin{equation*}
     (\triangle) ^{s}u_{1}(x) + cu_{1}'(x) + u_{1}(x) + \frac{\beta(x)}{u_{1}(x)}=t_{1} <\tilde{t}
\end{equation*}
hence $u_{1}$ is a strict $2\pi$-periodic supersolution of $(\ref{aps1})$ for $\tilde{t}$. By Lemma \ref{infcota} we have for $t\in[t_1,t_2]$
\begin{equation*}
   \left(\frac{\beta_{min}}{t}\right)^{1/\mu} \leq u(x), \qquad x\in \mathbb{R}.
\end{equation*}
In particular,
\begin{equation*}
   \left(\frac{\beta_{min}}{t_2}\right)^{1/\mu} \leq u(x), \qquad x\in \mathbb{R}.
\end{equation*}
Doing similarly as in the  proof of Lemma \ref{nosol} we get $r_{t_2}$ is a subsolution for $(\ref{aps1})$ for $\tilde{t}$. Therefore, reasoning like in step $3$ of the proof of the Theorem \ref{princ2} for $t=\tilde{t}$ we get
\begin{equation}\label{deg1}
   deg( I - KN, \Omega_{1}, 0)= 1,
\end{equation}
where $$\Omega_{1}= \{ u \in C^{1}_{2\pi}(\mathbb{R}): r_{t_2} < u(x) < u_{1}(x),\;  x\in\mathbb{R}\quad \|u'\|_{\infty}<N\}.$$ Then taking $R_{t_2}$ large enough that guarantee that $\Omega_1 \subseteq \Omega$, and by using the excision property, we conclude for $t=\tilde{t}$ that
\begin{equation*}
  deg( I -KN, \Omega\backslash \overline{\Omega_{1}}, 0) = deg( I-KN, \Omega, 0) - deg( I- KN, \Omega_{1}, 0) =-1,
\end{equation*}
that implies the existence of a $2\pi$-periodic solution of equation $(\ref{aps1})$ contained in $\Omega\backslash\overline{\Omega_{1}}$ for $\tilde{t}$. Again, since $\tilde{t}$ and $t_{2}$ are arbitrary, the proof is complete.
\fin
\end{desing}
\subsection{The Fractional Laplacian problem with attractive-repulsive singular nonlinearities}
In this subsection we provide some extension of the bounds found in the previous subsection, and we consider an equation with two singularities of the form,
 \begin{equation}\label{12}
  (\Delta)^s u(x) + cu'(x) + e(x) - \frac{\gamma (x)}{u^{\mu}(x)} + \frac{\beta (x)}{u^{\rho}(x)}= 0, \qquad x\in\mathbb{R},
\end{equation}
where $s\in(0,1)$, $c>0$, $\gamma$, $\beta$ and $e$ are Lipschitz continuous and $2\pi$-periodic, $\mu$ and $\rho$ are positive constants with $\mu \geq \rho$ and $\mu \geq 1$. We take advantage to prove the existence of positive solutions of $(\ref{12})$ using degree theory,  inspired by some recent works when $s=1$, (for example \cite{Lu}). Notice that $\gamma$ and $\beta$ can change signs which makes the proofs more interesting.\\
 From now on, we will use denotation $$\omega_{+}(x)=\max\{\omega(x),0\}, \qquad \omega_{-}(x)=-\min\{\omega(x),0\}.$$
Clearly,
\begin{equation*}
  \omega(x)=\omega_{+}(x) - \omega_{-}(x) \quad x\in \mathbb{R},\qquad \overline{\omega}=\overline{\omega_{+}} - \overline{\omega_{-}}.
\end{equation*}

\begin{teo}\label{e2sing}
Assume that $\bar{\gamma}>0$, $\bar{e}>0$ and that following condition holds
\begin{equation}\label{lim}
  \lim\limits_{x \to +0} \frac{cx}{2} - \frac{2\pi\overline{\gamma_{-}}}{x^{\mu}} - \frac{2\pi\overline{\beta_{+}}}{x^{\rho}} = + \infty,
\end{equation}
where $c>0$ is the that appears in ($\ref{12})$. Then, equation ($\ref{12}$) has at least one $2\pi$-periodic positive classical solution.
\end{teo}

To prove this Theorem we will use degree theory with some similar ideas to those of \cite{Lu} adapted to the nonlocal case. Consider now the family of equations
\begin{equation}\label{sing}
  (\Delta)^s u(x) + cu'(x) + \lambda e(x) - \lambda\frac{\gamma (x)}{u^{\mu}(x)} + \lambda\frac{\beta (x)}{u^{\rho}(x)}= 0 \quad \lambda\in (0,1].
\end{equation}
We need to find a priori estimates, for this, we will use the following Lemma that proves upper and lowers bound at some point.

\begin{lem}\label{casic}
Assume $\bar{\gamma}> 0$ and $\bar{e}>0$, then for each possible $2\pi$-periodic solution $u$ of (\ref{sing}) there exist $x_{0}, x_{1} \in [0,2\pi]$ such that
\begin{equation}\label{A0}
  u(x_{0}) \leq \max \left\{1,\left(\dfrac{\overline{\gamma_{+}} + \overline{\beta_{-}}}{\bar{e}}\right)^{\frac{1}{\rho}}\right\} = A_{0}
\end{equation}
and
\begin{equation}\label{A1}
  u(x_{1}) \geq A_{1},
\end{equation}
where
\begin{equation*}
  A_{1} = \begin{cases}
 \min\left\{\left(\dfrac{\overline{\gamma}}{2\overline{e_{+}}}\right)^{\frac{1}{\mu}}, \left(\dfrac{\overline{\gamma}}{2\overline{\beta_{+}}}\right)^{\frac{1}{\mu - \rho}}\right\}, & \text{si $\overline{\beta_{+}}>0$},\\
\left(\dfrac{\bar{\gamma}}{\overline{e_{+}}}\right)^{\frac{1}{\mu}}, & \text{si $\overline{\beta_{+}}=0$}.
\end{cases}
\end{equation*}
\end{lem}
\begin{dem}
  Let $u$ is a $2\pi$-periodic solution of ($\ref{sing}$). Then integrating the equation ($\ref{sing}$) in $[0,2\pi]$, we obtain
  \begin{equation*}
    \int\limits_{0}^{2\pi}(\Delta)^s u(x) + c \int\limits_{0}^{2\pi}u'(x) + \lambda\int\limits_{0}^{2\pi} e(x) - \lambda\int\limits_{0}^{2\pi}\frac{\gamma (x)}{u^{\mu}(x)} + \lambda\int\limits_{0}^{2\pi}\frac{\beta (x)}{u^{\rho}(x)}= 0
  \end{equation*}

  By \cite[Lemma 2.2]{lis}, together with the $2\pi$-periodicity of $u$, we get
  \begin{eqnarray}\label{cot}
  0< 2\pi\bar{e}&=& \lambda\int\limits_{0}^{2\pi}\frac{\gamma (x)}{u^{\mu}(x)} - \lambda\int\limits_{0}^{2\pi}\frac{\beta (x)}{u^{\rho}(x)} \nonumber \\
       &\leq& \int\limits_{0}^{2\pi}\frac{\gamma_{+} (x)}{u^{\mu}(x)} + \int\limits_{0}^{2\pi}\frac{\beta_{-} (x)}{u^{\rho}(x)}.
\end{eqnarray}
  From the above, we can conclude that there is a point $x_{0}\in [0,2\pi]$ such that $(\ref{A0})$ is true. Indeed, if ($\ref{A0}$) does not hold, then

  \begin{equation}\label{4}
    u(x)> 1  \quad \text{for all} \; x\in [0,2\pi],
  \end{equation}
 and
 \begin{equation}\label{5}
   u(x)> \left(\dfrac{\overline{\gamma_{+}} + \overline{\beta_{-}}}{\bar{e}}\right)^{\frac{1}{\rho}} \quad \text{for all} \; x\in [0,2\pi].
 \end{equation}

  Since $\mu\geq \rho$ and by ($\ref{4}$), using the Mean Value Theorem of integrals in ($\ref{cot}$), we have that there exists a point $x_{0}\in [0,2\pi]$ such that
 \begin{equation*}
   2\pi \bar{e} \leq \frac{2 \pi(\overline{\gamma_{+}} +\overline{\beta_{-}})}{u^{\rho}(x_{0})}
 \end{equation*}
which implies that
\begin{equation*}
  u(x_{0}) \leq \left(\dfrac{\overline{\gamma_{+}} + \overline{\beta_{-}}}{\bar{e}}\right)^{\frac{1}{\rho}},
\end{equation*}
 this contradicts ($\ref{5}$), so ($\ref{A0}$) holds.\\

 Now, we shall proof that ($\ref{A1}$) is true. Indeed, multiplying ($\ref{sing}$) for $u^{\mu}$ and integrating the equation ($\ref{sing}$) in $[0,2\pi]$ we get
  \begin{equation*}
    \int\limits_{0}^{2\pi}(\Delta)^s u(x)u^{\mu}(x)dx + c \int\limits_{0}^{2\pi}u'(x)u^{\mu}(x)dx + \lambda\int\limits_{0}^{2\pi} e(x)u^{\mu}(x)dx - \lambda\int\limits_{0}^{2\pi}\gamma(x)dx + \lambda\int\limits_{0}^{2\pi} \beta (x)u^{\mu - \rho}(x)dx= 0,
  \end{equation*}
that is,
 \begin{equation}\label{6}
  \lambda 2\pi\bar{\gamma} +  \int\limits_{0}^{2\pi}(-\Delta)^s u(x)u^{\mu}(x)dx - c \int\limits_{0}^{2\pi}u'(x)u^{\mu}(x)dx= \lambda\int\limits_{0}^{2\pi} \beta (x)u^{\mu - \rho}(x)dx + \lambda\int\limits_{0}^{2\pi} e(x)u^{\mu}(x)dx.
  \end{equation}
Since,

\begin{eqnarray}\label{7}
  \int\limits_{0}^{2\pi}(-\Delta)^s u(x)u^{\mu}(x)dx &=& \int\limits_{0}^{2\pi}(-\Delta)^s u(x)u(x)u^{\mu-1}(x)dx\nonumber\\
                                                    &\geq& u_{min}^{\mu - 1}[u]_{H^{s}(0,2\pi)}^{2} \geq 0.
\end{eqnarray}
Using the fact that
\begin{equation*}
  - c \int\limits_{0}^{2\pi}u'(x)u^{\mu}(x) = 0,
\end{equation*}
 we have
 \begin{eqnarray}\label{8}
  2\pi \bar{\alpha}  &\leq& \int\limits_{0}^{2\pi} \beta (x)u^{\mu - \rho}(x) + \int\limits_{0}^{2\pi} e(t)u^{\mu}(x)\nonumber\\
                                                   &\leq& \int\limits_{0}^{2\pi} \beta_{+}(x)u^{\mu - \rho}(x) +  \int\limits_{0}^{2\pi}e_{+}(x)u^{\mu}(x)
\end{eqnarray}

 Thus, reasoning as before we can conclude that there is a point $x_{1}\in [0,2\pi]$ such that
 \begin{equation*}
   u(x_{1}) \geq A_{1}.
 \end{equation*}
 \fin
\end{dem}\\

\begin{lem}\label{cotainf}
Let us assume that $\bar{\gamma}>0$, $\bar{e}>0$ and $(\ref{lim})$ hold. Then there exists $r>0$ such that each $2\pi$-periodic solution of ($\ref{sing}$) satisfies $u(t) > r$ for all $x\in[0,2\pi]$.
\end{lem}

\begin{dem}
Let $x^{n}_{3}$ be the minimum point of $u_{n}$ in $[0,2\pi]$. We assume by contradiction that $u_{n}(x^{n}_{3}) \to 0$ as $n\to \infty$.\\
 Then, we suppose that  $x^{n}_{3}< x_{0}$, where $x_{0}$ is given in Lemma \ref{casic}, and integrating ($\ref{sing}$) over $[x^{n}_{3}, x_{0}]$, we obtain
 \begin{equation}\label{inter1}
    \int\limits_{x^{n}_{3}}^{x_{0}}(\Delta)^s u_{n}(x)dx + c \int\limits_{x^{n}_{3}}^{x_{0}}u_{n}'(x)dx= \lambda\int\limits_{x^{n}_{3}}^{x_{0}}\frac{\gamma (x)}{u_{n}^{\mu}(x)}dx - \lambda\int\limits_{x^{n}_{3}}^{x_{0}}\frac{\beta (x)}{u_{n}^{\rho}(x)}dx - \lambda\int\limits_{x^{n}_{3}}^{x_{0}}e(x)dx.
 \end{equation}
 Since
 \begin{equation*}
  \int\limits_{x^{n}_{3}}^{x_{0}}u_{n}'(x)dt=u_{n}(x_{0}) -u_{n}(x^{n}_{3}) \geq 0,
 \end{equation*}
then
\begin{eqnarray}\label{lapla}
   \int\limits_{x^{n}_{3}}^{x_{0}}(\Delta)^s u_{n}(x)dx &\leq&  \lambda\int\limits_{x^{n}_{3}}^{x_{0}}\frac{\gamma (x)}{u_{n}^{\mu}(x)}dx - \lambda\int\limits_{x^{n}_{3}}^{x_{0}}\frac{\beta (x)}{u_{n}^{\rho}(x)}dx -\lambda \int\limits_{x^{n}_{3}}^{x_{0}}e(x)dx \nonumber\\
   &\leq&  \lambda\int\limits_{0}^{2\pi}\frac{\gamma_{+} (x)}{u_{n}^{\mu}(x)}dx + \lambda\int\limits_{0}^{2\pi}\frac{\beta_{-} (x)}{u_{n}^{\rho}(x)}dx +\lambda\int\limits_{0}^{2\pi}e_{-}(x)dx.
\end{eqnarray}
It is important to remark that, if $x^{n}_{3}> x_{0}$, by the $2\pi$-periodicity of $u$ we can take $x^{n}_{*}$ be the  minimum point in $[-2\pi,0]$ and integrating over $[x^{n}_{*},x_{0}]$.\\
On the other hand, integrating ($\ref{sing}$) on $[0,2\pi]$ and using \cite[Lemma 2.2]{lis} and the $2\pi$-periodicity of $u'$, we get
\begin{equation*}
\lambda\int\limits_{0}^{2\pi}\frac{\gamma (x)}{u^{\mu}(x)}dx =\lambda\int\limits_{0}^{2\pi} \frac{\beta (x)}{u^{\rho}(x)}dx +\lambda \int\limits_{0}^{2\pi} e(x)dx,
\end{equation*}
that is,
\begin{equation}\label{m}
 \lambda \int\limits_{0}^{2\pi}\frac{\gamma_{+}(x)}{u^{\mu}(x)}dx + \lambda\int\limits_{0}^{2\pi} \frac{\beta_{-}(x)}{u^{\rho}(x)}dx + \lambda\int\limits_{0}^{2\pi} e_{-}(x)dx= \lambda \int\limits_{0}^{2\pi}\frac{\gamma_{-}(x)}{u^{\mu}(x)}dx + \lambda \int\limits_{0}^{2\pi} \frac{\beta_{+}(x)}{u^{\rho}(x)}dx + \lambda\int\limits_{0}^{2\pi} e_{+}(x)dx.
\end{equation}
Substituting it into $(\ref{lapla})$, we have

\begin{equation}\label{d}
   \int\limits_{x^{n}_{3}}^{x_{0}}(\Delta)^s u_{n}(x) \leq \frac{\lambda 2\pi\overline{\gamma_{-}}}{u_{n}^{\mu}(x_{3}^{n})} + \frac{ \lambda 2\pi\overline{\beta_{+}}}{u_{n}^{\rho}(x_{3}^{n})} +  \lambda 2\pi\overline{e_{+}}.
\end{equation}

Hence, from $(\ref{inter1})$ together with ($\ref{d}$) gives
\begin{equation*}
  c\int\limits_{x^{n}_{3}}^{x_{0}}u_{n}'(x)dx \geq \lambda \int\limits_{x^{n}_{3}}^{x_{0}}\frac{\gamma (x)}{u_{n}^{\mu}(x)}dx - \lambda\int\limits_{x^{n}_{3}}^{x_{0}} \frac{\beta (x)}{u_{n}^{\rho}(x)}dx - \lambda\int\limits_{x^{n}_{3}}^{x_{0}} e(x)dx - \frac{\lambda 2\pi\overline{\gamma_{-}}}{u_{n}^{\mu}(x_{3}^{n})} - \frac{\lambda 2\pi\overline{\beta_{+}}}{u_{n}^{\rho}(x_{3}^{n})} - \lambda 2\pi\overline{e_{+}} ,
\end{equation*}
so that
\begin{eqnarray*}
  cu_{n}(x_{0}) &\geq& cu_{n}(x^{n}_{3}) + \lambda\int\limits_{x^{n}_{3}}^{x_{0}}\frac{\gamma (x)}{u_{n}^{\mu}(x)}dx -\lambda\int\limits_{x^{n}_{3}}^{x_{0}} \frac{\beta (x)}{u_{n}^{\rho}(x)}dx - \lambda\int\limits_{x^{n}_{3}}^{x_{0}} e(x) dx\nonumber \\
  && \qquad\qquad\qquad\qquad\qquad\qquad- \frac{\lambda 2\pi\overline{\gamma_{-}}}{u_{n}^{\mu}(x_{3}^{n})} - \frac{\lambda 2\pi\overline{\beta_{+}}}{u_{n}^{\rho}(x_{3}^{n})} - \lambda 2\pi\overline{e_{+}}\\
  &\geq&  \lambda cu_{n}(x^{n}_{3}) - \lambda\int\limits_{0}^{2\pi}\frac{\gamma_{-}(x)}{u_{n}^{\mu}(x)}dx -\lambda\int\limits_{0}^{2\pi} \frac{\beta_{+} (x)}{u_{n}^{\rho}(x)}dx - \lambda\int\limits_{0}^{2\pi} e_{+}(x)dx\nonumber \\
   && \qquad\qquad\qquad\qquad\qquad\qquad - \frac{\lambda 2\pi\overline{\gamma_{-}}}{u_{n}^{\mu}(x_{3}^{n})} - \frac{\lambda 2\pi\overline{\beta_{+}}}{u_{n}^{\rho}(x_{3}^{n})} - \lambda 2\pi\overline{e_{+}}\\
   &\geq& 2\lambda\left(\frac{c}{2}u_{n}(x^{n}_{3})- \frac{2\pi\overline{\gamma_{-}}}{u_{n}^{\mu}(x_{3}^{n})} - \frac{2\pi\overline{\beta_{+}}}{u_{n}^{\rho}(x_{3}^{n})} - 2\pi\overline{e_{+}}\right).
\end{eqnarray*}
 By ($\ref{lim}$) we have $u_{n}(x_{0}) \to +\infty$, but this contradicts (\ref{A0}). Hence, the result follows.
\fin
\end{dem}\\

\begin{lem}\label{cotasup}
Let us assume that $\bar{\gamma}>0$ and $u$ a solution positive $2\pi$-periodic of ($\ref{sing}$) satisfying $(\ref{lim})$, then there exists a constant $R>0$ such that
\begin{equation*}
  0< u(x) <R \qquad x\in[0,2\pi].
\end{equation*}
\end{lem}
\begin{dem}
Let $x_{3}$ be defined as the same as the one in the proof of Lemma \ref{cotainf}, that is $u(x_{3})=\min\limits_{x\in[0,2\pi]}u(x)$. By Lemma \ref{casic}, it is clear that
\begin{eqnarray}\label{l2}
  u(x) &\leq& u(x_{0}) + \int\limits_{0}^{2\pi}|u'| \nonumber \\
   &\leq& A_{0} + \sqrt{2\pi}\|u'\|_{L^{2}(0,2\pi)}.
\end{eqnarray}
Now, multiplying ($\ref{sing}$) for $u'$ and integrating the equation
\begin{equation*}
    \int\limits_{0}^{2\pi}(\Delta)^s u(x)u'(x)dx + c \int\limits_{0}^{2\pi}|u'(x)|^{2}dx=  \lambda\int\limits_{0}^{2\pi}\frac{\gamma(x)}{u^{\mu}(x)}u'(x)dx - \lambda\int\limits_{0}^{2\pi} \frac{\beta (x)}{u^{\rho}(x)}u'(x)dx- \lambda\int\limits_{0}^{2\pi} e(x)u'(x)dx.
  \end{equation*}
  By \cite[Lemma 4.2]{lis} and using H$\ddot{o}$lder inequality together with the Lemma \ref{cotainf} and $0<\lambda\leq 1$ gives
  \begin{eqnarray*}
     c\|u'\|_{L^{2}(0,2\pi)}^{2} &\leq&\int\limits_{0}^{2\pi}\frac{|\gamma|}{u^{\mu}(x_{3})}|u'| + \int\limits_{0}^{2\pi} \frac{|\beta|}{u^{\rho}(x_{3})}|u'|+ \int\limits_{0}^{2\pi} |e| |u'| \\
     &\leq& C\|u'\|_{L^{2}(0,2\pi)} (\|\gamma\|_{L^{2}(0,2\pi)} + \|\beta\|_{L^{2}(0,2\pi)} + \|e\|_{L^{2}(0,2\pi)}),
  \end{eqnarray*}
  so that,

\begin{equation*}
  \|u'\|_{L^{2}(0,2\pi)}\leq C(\|\gamma\|_{L^{2}(0,2\pi)} + \|\beta\|_{L^{2}(0,2\pi)} + \|e\|_{L^{2}(0,2\pi)}),
\end{equation*}
  which, together with ($\ref{l2}$) implies
  \begin{equation*}
    u(x) \leq A_{0} + C(\|\gamma\|_{L^{2}(0,2\pi)} + \|\beta\|_{L^{2}(0,2\pi)} + \|e\|_{L^{2}(0,2\pi)}):=R.
  \end{equation*}
\fin
\end{dem}\\

\begin{obs}\label{obs}
  If $\bar{\gamma}>0$ and $\bar{e}>0$, then there are two constants $M_{1}$ and $M_{2}$ with $0<M_{1}<M_{2}$ such that
  \begin{equation*}
     \frac{\bar{\gamma}}{x^{\mu}} -\frac{\bar{\beta}}{x^{\rho}} - \bar{e}>0, \quad \text{for all} \; x\in(0,M_{1}),
  \end{equation*}
  and
   \begin{equation*}
     \frac{\bar{\gamma}}{x^{\mu}} -\frac{\bar{\beta}}{x^{\rho}} - \bar{e}<0, \quad \text{for all} \; x\in(M_{2}, + \infty)
  \end{equation*}
\end{obs}

We now proof the existence of solutions for $(\ref{12})$. For this, we will use the ideas of the Continuation Theorem of \cite[Theorem 2.4]{Ma}.

\begin{demstwo}
Let $r_{0}= \min\{M_{1},r\}$ and $R_{0}= \max\{M_{2},R\}$, where $M_{1}$ and $M_{2}$ were determined in Remark \ref{obs}, we have
\begin{equation}\label{ct}
     \frac{\bar{\gamma}}{a^{\mu}} -\frac{\bar{\beta}}{a^{\rho}} - \bar{e}>0, \quad \text{for all} \; a\in(0,r_{0}]
  \end{equation}
  and
   \begin{equation}\label{ctt}
     \frac{\bar{\gamma}}{a^{\mu}} -\frac{\bar{\beta}}{a^{\rho}} - \bar{e}<0, \quad \text{for all} \; a\in [R_{0}, + \infty),
  \end{equation}
hence each possible constant solution of the equation
\begin{equation}\label{csolu}
   \frac{\bar{\gamma}}{\tilde{a}^{\mu}} -\frac{\bar{\beta}}{\tilde{a}^{\rho}} - \bar{e}=0.
\end{equation}
 satisfies $r_{0}<\tilde{a}<R_{0}$.
We can define a map $K: L^{\infty}(\mathbb{R}) \to C^{1,\alpha}_{2\pi}(\mathbb{R})$ by $K(z)= u$ where $u$ is a solution of
\begin{equation*}
  (\triangle)^{s}u(x) + cu'(x)= z(x),
\end{equation*}
where $K$ is compact.\\

Let us now define the map $N: C^{1}_{2\pi}(\mathbb{R}) \to L^{\infty}(\mathbb{R})$ by
\begin{equation}\label{defnn}
  Nu= \frac{\gamma(\cdot)}{u^{\mu}(\cdot)} - \frac{\beta(\cdot)}{u^{\rho}(\cdot)} -  e(\cdot).
\end{equation}
Reasoning like in the proof of Lemma \ref{cotac1al} we have that there exists $N=N(R,\gamma,\beta,e,c)$ and $\alpha\in (0,1)$ such that
\begin{equation*}
  \|u\|_{C^{1,\alpha}(\mathbb{R})}< N.
\end{equation*}
Let $\Omega= \{ u\in C^{1}_{2\pi}(\mathbb{R}) : r_{0} <u(x) <R_{0},\;  x\in [0,2\pi] \;  \|u'\|_{\infty}< N\}$ and define the continuous projectors $Q:  L^{\infty}(\mathbb{R}) \to L^{\infty}(\mathbb{R})$ by the constant function
 \begin{equation*}
   Qy= \frac{1}{2\pi}   \int\limits_{0}^{2\pi} y(t)dt.
 \end{equation*}

Let us define one parameter family of problems
\begin{equation*}
   u=K((1-\lambda)QNu + \lambda Nu),  \;\; \lambda \in [0,1].
\end{equation*}
Explicitly,
\begin{equation}\label{homo}
   (\triangle)^{s}u(x) + cu'(x)= (1-\lambda)QNu + \lambda Nu.
\end{equation}
For $\lambda \in [0, 1]$, applying $Q$ to both members of $(\ref{homo})$ we have by \cite[Lemma 2.2]{lis} and the $2\pi$-periodicity of $u'$
\begin{equation}\label{fin}
  \frac{1}{2\pi}\int\limits_{0}^{2\pi} \left(\frac{\gamma (x)}{u^{\mu}(x)} - \frac{\beta (x)}{u^{\rho}(x)} - e(x)\right) dx=0.
\end{equation}
Since also we obtain $(\ref{fin})$ when applying $Q$ in $(\ref{sing})$, then for all $\lambda\in (0, 1]$, problem $(\ref{sing})$ and problem $(\ref{homo})$ are equivalents. Hence, Lemmas \ref{cotainf} and \ref{cotasup} imply that $(\ref{homo})$ does not have a solution on $\partial\Omega \times (0,1]$. For $\lambda=0$, $(\ref{homo})$ is equivalent to the problem
\begin{equation}\label{qn}
   (\triangle)^{s}u(x) + cu'(x)=  \frac{1}{2\pi}\int\limits_{0}^{2\pi}\left(\frac{\gamma (x)}{u^{\mu}(x)} - \frac{\beta (t)}{u^{\rho}(x)} - e(x) \right)dx,
\end{equation}
then, applying $Q$ to both members of this equation and by \cite[Lemma 2.2]{lis} we get
\begin{equation*}
  QNu=0 \qquad  (\triangle)^{s}u(x) + cu'(x)= 0,
\end{equation*}
  reasoning similar to the proof of \cite[Theorem 1.4]{lis} we obtain that $u$ is constant, since the beginning of the proof, it was obtained that the constant solutions of $QN\tilde{a}=0$ satisfy the inequality $r_{0}<\tilde{a}<R_{0}$ we have that $(\ref{homo})$ has no solution on $\partial\Omega\times[0,1]$.  Therefore, the $deg(I- K((1-\lambda)QN + \lambda N), \Omega, 0)$  is well defined  for all $\lambda\in[0,1]$ and by homotopy invariant of the degree we have
 \begin{equation*}
   deg(I - KN, \Omega, 0)= deg(I - KQN, \Omega, 0)=deg(I - KQN, \Omega\cap E, 0),
 \end{equation*}
where $E\subset C^{1,\alpha}_{2\pi}(\mathbb{R}):$ is the one dimensional space of constants maps, and the last equality is since all solutions of $(\ref{qn})$ are constants maps as proved above. From here, by \eqref{ct} and \eqref{ctt} we get $deg(I - KQN, \Omega\cap E, 0)\neq 0$. Thus, we can conclude  that the equation $(\ref{12})$ has at least one positive classical solution $2\pi$-periodic.
 \fin
\end{demstwo}\\

{\bf Acknowledgements}
B. B. was partially supported by AEI Grant MTM$2016$- $80474$-P and Ramón y Cajal fellowship RYC$2018$-$026098$-I (Spain). L. C. was partially supported by ANID \# 21191475. A. Q. was partially supported by FONDECYT Grant \# 1190282.
\bibliographystyle{amsplain}
\bibliography{periodicAPS}

\providecommand{\bysame}{\leavevmode\hbox to3em{\hrulefill}\thinspace}
\providecommand{\MR}{\relax\ifhmode\unskip\space\fi MR }
\providecommand{\MRhref}[2]{%
  \href{http://www.ams.org/mathscinet-getitem?mr=#1}{#2}
}
\providecommand{\href}[2]{#2}
\begin{thebibliography}{10}

\bibitem{ambrosetti}
A.~Ambrosetti and G.~Prodi, \emph{On the inversion of some differentiable
  mappings with singularities between banach spaces}, Annali di Matematica Pura
  ed Applicata \textbf{93} (1972), 231--246.

\bibitem{ambrosioV}
V.~Ambrosio, \emph{An ambrosetti–prodi type result for fractional spectral
  problems}, Mathematische Nachrichten \textbf{293} (2020).

\bibitem{ambrosioymoli}
V.~Ambrosio and G.~Molica~Bisci, \emph{Periodic solutions for a fractional
  asymptotically linear problem}, Proceedings of the Royal Society of
  Edinburgh: Section A Mathematics \textbf{149} (2019), no.~3, 593–615.

\bibitem{quaas}
B.~Barrios, J.~García-Melián, and A.~Quaas, \emph{Periodic solutions for the
  one-dimensional fractional laplacian}, Journal of Differential Equations
  \textbf{267} (2019), no.~9, 5258--5289.

\bibitem{berger}
M.S. Berger and E.~Podolak, \emph{On the solutions of a nonlinear dirichlet
  problem}, Indiana Univ. Math. J. \textbf{24} (1975), 837--846.

\bibitem{biswas}
A.~Biswas and J.~L\"orinczi, \emph{Ambrosetti–prodi type results for
  dirichlet problems of fractional laplacian-like operators}, Integral
  Equations and Operator Theory \textbf{92} (2020), no.~26.

\bibitem{caffarelli}
L.~Caffarelli and L.~Silvestre, \emph{Regularity theory for fully nonlinear
  integro-differential equations}, Communications on Pure and Applied
  Mathematics \textbf{62} (2009).

\bibitem{Luisxavier}
L.~A. Caffarelli and X.~Cabre, \emph{Fully nonlinear elliptic equations,
  american mathematical society colloquium publications.}, American
  Mathematical Society Colloquium Publications. \textbf{43} (1995).

\bibitem{lis}
L.~Carrero and A.~Quaas, \emph{Periodic solutions for one-dimensional nonlinear
  nonlocal problem with drift including singular nonlinearities}, Proceedings
  of the Royal Society of Edinburgh: Section A Mathematics, 1--33.

\bibitem{chia}
R.~Chiappinelli, J.~Mawhin, and R.~Nugari, \emph{Generalized ambrosetti-prodi
  conditions for nonlinear two-point boundary value problems}, Journal of
  Differential Equations \textbf{69} (1987), no.~3, 422--434.

\bibitem{zcg}
Z.~Du and C.~Gui, \emph{Further study on periodic solutions of elliptic
  equations with a fractional laplacian},  \textbf{193} (2020), 111417.

\bibitem{davilaqt}
G.~Dávila, A.~Quaas, and E.~Topp, \emph{Existence, nonexistence and
  multiplicity results for nonlocal dirichlet problems}, Journal of
  Differential Equations \textbf{266} (2018), no.~9, 5971-5997.

\bibitem{mawhin}
C.~Fabry, J.~Mawhin, and M.~Nkashama, \emph{A multiplicity result for periodic
  solutions of forced nonlinear second order ordinary differential equations},
  Bulletin of the London Mathematical Society \textbf{18} (1986), no.~2,
  173--180.

\bibitem{felmer}
P.~Felmer and A.~Quaas, \emph{Fundamental solutions and liouville type theorems
  for nonlinear integral operators}, Advances in Mathematics \textbf{226}
  (2011), no.~3, 2712--2738.

\bibitem{fernandez}
X.~Fernández-Real and X.~Ros-Oton, \emph{The obstacle problem for the
  fractional laplacian with critical drift}, Mathematische Annalen \textbf{371}
  (2017).

\bibitem{Ma}
P.~Fitzpatrick, M.~Martelli, J.~Mawhin, and R.~Nussbaum, \emph{Lecture on
  topological methods for ordinary differential equations}, Springer-Verlag
  (1993), 218.

\bibitem{at}
A.~Gutiérrez and P.~Torres, \emph{Non-autonomous saddle-node bifurcation in a
  canonical electrostatic mems}, International Journal of Bifurcation and Chaos
  \textbf{23} (2013), no.~05, 350088.

\bibitem{Nikos}
N.~Katzourakis, \emph{An introduction to viscosity solutions for fully
  nonlinear 2nd order pde and applications to calculus of variations in
  $l^\infty$}, Lecture Notes in Mathematics -Springer-verlag- (2014).

\bibitem{interpol}
N.~Krylov, \emph{Lectures on elliptic and parabolic equations in holder
  spaces}, 1996.

\bibitem{Laskin}
N.~Laskin, \emph{Fractional quantum mechanics and l{\'e}vy path integrals},
  Physics Letters A \textbf{268} (1999), 298--305.

\bibitem{Lu}
S.~Lu and X.~Yu, \emph{Periodic solutions for second order differential
  equations with indefinite singularities}, Advances in Nonlinear Analysis
  \textbf{9} (2019), 994 -- 1007.

\bibitem{principal}
J.~Mawhin, \emph{The periodic ambrosetti-prodi problem for nonlinear
  perturbations of the p-laplacian}, Journal of The European Mathematical
  Society - J EUR MATH SOC \textbf{8} (2006), 375--388.

\bibitem{Ralf}
R.~Metzler and J.~Klafter, \emph{The restaurant at the end of the random walk:
  recent developments in the description of anomalous transport by fractional
  dynamics},  \textbf{37} (2004), no.~31, R161--R208.

\bibitem{Lsilvestreshw}
R.~W. Schwab and L.~Silvestre, \emph{Regularity for parabolic
  integro-differential equations with very irregular kernels}, Analysis and PDE
  \textbf{9} (2014).

\bibitem{silvestre}
L.~Silvestre, \emph{Regularity of the obstacle problem for a fractional power
  of the laplace operator}, Communications on Pure and Applied Mathematics
  \textbf{60} (2007), 67 -- 112.

\bibitem{Lsilvestre}
\bysame, \emph{On the differentiability of the solution to an equation with
  drift and fractional diffusion}, Indiana University Mathematics Journal
  \textbf{61} (2010).

\bibitem{fisica2}
G.~M. Viswanathan, V.~Afanasyev, E.J. Buldyrev, S.V.~Murphy, P.A. Prince, and
  H.E. Stanley, \emph{Lévy flight search patterns of wandering albatrosses},
  Nature \textbf{318} (1996), 413--415.

\bibitem{fisica}
G.~Zhang and B.~Li, \emph{Anomalous vibrational energy diffusion in carbon
  nanotubes.}, The Journal of chemical physics \textbf{123 1} (2005), 014705.

\bibitem{zhuo}
R.~Zhuo, W.~Chen, X.~Cui, and Z.~Yuan, \emph{A liouville theorem for the
  fractional laplacian}, arXiv e-prints (2014).

\end{thebibliography}
\end{document}